\documentclass[11pt,reqno]{amsart}
\usepackage[T1]{fontenc}
\usepackage{lmodern,amsmath,amssymb,amsthm,geometry,microtype,needspace}
\usepackage[colorlinks=true,linkcolor=blue,citecolor=blue,urlcolor=blue]{hyperref}
\hypersetup{pdftitle={One-fibre determination and maximal domains for projectively flat Finsler metrics},pdfauthor={Benling Li, Kefeng Liu, Liangdi Zhang}}
\makeatletter
\@ifundefined{subjclassname@2020}{%
  \@namedef{subjclassname@2020}{\textup{2020} Mathematics Subject Classification}%
}{}
\makeatother
\numberwithin{equation}{section}
\newtheorem{theorem}{Theorem}[section]
\newtheorem{lemma}[theorem]{Lemma}
\newtheorem{proposition}[theorem]{Proposition}
\newtheorem{corollary}[theorem]{Corollary}
\newtheorem{theoremA}{Theorem}

\newtheorem{theoremB}{Theorem}

\newtheorem{theoremC}{Theorem}

\theoremstyle{definition}

\theoremstyle{remark}
\newtheorem{remark}[theorem]{Remark}

\newcommand{\R}{\mathbb R}\newcommand{\C}{\mathbb C}
\newcommand{\dd}{\mathrm d}\newcommand{\Hess}{\operatorname{Hess}}
\newcommand{\Id}{\operatorname{Id}}\newcommand{\im}{\operatorname{Im}}
\newcommand{\re}{\operatorname{Re}}\newcommand{\calO}{\mathcal O}
\newcommand{\calH}{\mathcal H}\newcommand{\Sn}{\mathbb S^{n-1}}

\newcommand{\Dop}{\mathcal D}
\title[One fibre and maximal domains]{One-fibre determination and
maximal domains for projectively flat Finsler metrics}
\author{Benling Li}
\address{School of Mathematics and Statistics, Ningbo University,
315211 Ningbo, China}
\email{libenling@nbu.edu.cn}
\author{Kefeng Liu}
\address{Mathematical Sciences Research Center, Chongqing University of
Technology, Chongqing 400054, China}
\email{kefengliu@cqut.edu.cn}
\author{Liangdi Zhang}
\address{Mathematical Sciences Research Center, Chongqing University of
Technology, Chongqing 400054, China}
\email{ldzhang91@163.com}
\date{}
\begin{document}

\begin{abstract}
In the regular case, Hilbert's fourth problem asks for the classification
of projectively flat Finsler metrics. For constant flag curvature, the local classification determines the metric from the restrictions of
the metric and its projective factor to one tangent fibre; in nonpositive curvature this is unconditional, while in positive
curvature it requires the metric to admit a complex extension. We remove this requirement: smooth solutions are automatically
real analytic, and smooth data admit a local solution exactly when real
analytic. Hence the metric and its projective factor on one fibre determine every solution on any connected
coordinate domain, with no completeness assumption. Without completeness
we also determine the maximal connected base domain on which the equivalent scalar
equation has a smooth solution in every nonzero direction; it is convex,
explicit, and contains the strong-convexity locus.
For equal Randers data we compute this locus in every dimension: at and
beyond a critical parameter it consists of two convex components.
For general data, however, a component of the locus need not be convex
in dimensions at least three; we give a Hilbert-metric example.
\end{abstract}

\maketitle

\section{Introduction}
Hilbert's fourth problem asks for geometries in which straight segments are shortest paths~\cite[Problem~4]{Hilbert1902}.
In the regular case,
this leads to characterizing the projectively flat Finsler metrics,
that is, the metrics whose geodesics are straight lines.
The projective viewpoint goes back to Beltrami, who showed that a surface metric can be mapped to the plane with geodesics becoming straight lines exactly when it has constant curvature~\cite{Beltrami1865}.
For $n\geq 2$, let $\Omega\subset\R^n$  be an open set
containing the origin, and let $F = F(x,y)$ be a Finsler metric on $\Omega$ that is
smooth on the slit bundle $\Omega\times(\R^n\setminus\{0\})$.
Hamel~\cite{Hamel1903} expressed the straight-geodesic condition as the
differential equations
\begin{equation}\label{eq:hamel}
 \Dop F_{y^k}=F_{x^k}, \qquad \Dop = y^j\partial_{x^j}.
\end{equation}
Then the geodesic coefficients take the form $G^i=Py^i$, where $P=P(x,y)$ is the  projective factor of $F$.
In the broader setting of metric spaces, the problem was developed by \mbox{Busemann}~\cite{Busemann1976}, Pogorelov~\cite{Pogorelov1979}, and Szab\'o~\cite{Szabo1986}; see Papadopoulos~\cite{Papadopoulos2014} for a survey.
Integral geometry gives a complementary description:
Alexander~\cite{Alexander1988} related Crofton representations to
zonoids, \'Alvarez Paiva and Fernandes~\cite{AlvarezFernandes1998}
extended Crofton formulas to projective Finsler spaces, and \'Alvarez
Paiva~\cite{AlvarezPaiva2005} developed the symplectic side of the
theory.

The constant-flag-curvature problem imposes a further differential constraint.
Hamel's equation together with constant flag curvature $\mathbf{K}$ gives
\begin{equation}\label{eq:system}
 F_{x^k}=(PF)_{y^k},\qquad P_{x^k}=PP_{y^k}-\mathbf{K}FF_{y^k}.
\end{equation}
These are the Berwald equations; they are equivalent to projective
flatness together with constant flag curvature (see Berwald~\cite{Berwald1929}).
Rescaling $F$ multiplies $\mathbf K$ by a positive constant, so it
suffices to consider $\mathbf K\in\{-1,0,1\}$.
Funk settled the two-dimensional case \cite{Funk1929}, and Berwald the
$n$-dimensional one \cite{Berwald1929}.
In 2003, Shen~\cite{Shen} gave implicit algebraic equations for projectively flat Finsler metrics of constant flag curvature, formulating the local problem in terms of the initial metric and projective factor.
In 2014, Li~\cite{Li} proved the converse, showing that the solution is
determined by its restriction to a single tangent fibre; for
$\mathbf{K}\le0$ this holds without further hypothesis.

Non-Riemannian examples exist for each value of the constant curvature: Berwald's metric and non-Euclidean Minkowski norms for $\mathbf{K}=0$, the Funk and Hilbert metrics for $\mathbf{K}<0$, and Bryant's metrics on $\mathbb{S}^2$ for $\mathbf{K}=1$ \cite{Bryant1996, Bryant1997}. Further examples and geodesic properties are discussed in \cite{Bryant, Bryant2006, ChengLi, BFIMZ}.

Building upon the local classification, Li and Zhao~\cite{LZ} systematically studied the global classification of these metrics. For \(\mathbf{K} \le 0\), they classified the metrics under the assumption of \emph{forward completeness}, which for negative curvature naturally restricts the base domain to be strictly convex. For \(\mathbf{K}=1\), they obtained classification results via \emph{metric completion}. These global works heavily rely on the completeness of the metric or the global convexity of the domain. In particular, they proved that in dimension two every connected
component of the regularity locus is convex, and described the degenerate
directions at its boundary; the analysis of the present paper starts from
these results.

In this paper, we drop the completeness assumption.
We focus on the local and global behaviour of the underlying partial
differential equations and on the geometric positivity of the metric.
Specifically, we address three questions: where the PDE system \eqref{eq:system} admits a unique
solution, what the maximal domain of existence for the scalar equations
is, and where the reconstructed metric remains strongly convex.

In positive curvature, the constructions of Shen~\cite{Shen} and
Li~\cite{Li} require the metric to admit a complex extension. Our
first result shows that smooth solutions are automatically real analytic,
and identifies exactly which smooth data admit a local solution.

For $\mathbf{K}=1$, set $Q=P+iF$. Then
\[
 Q_{x^k}=QQ_{y^k},\qquad \im Q>0.
\]
Near a point over $x_0\in\Omega$, the change of real coordinates \(z=y+(x-x_0)Q\) turns this system
into the Cauchy--Riemann equations. Section~\ref{sec:holo} proves this
step and uses the holomorphic implicit function theorem to recover
analyticity in \((x,y)\); see also
\cite{HormanderComplex,KrantzParksImplicit} for the analytic tools.

\begin{theoremA}\label{thmA}
Let $F$ be a smooth projectively flat Finsler metric of
flag curvature $\mathbf{K}=1$ on a connected open set $\Omega\subset\R^n$, with projective factor $P$. Then
$F$ and $P$ are real analytic on $\Omega\times(\R^n\setminus\{0\})$.
Conversely, let $\psi$ be a smooth Minkowski norm and let $\phi$ be a
smooth positively one-homogeneous function on $\R^n\setminus\{0\}$. The
data $F(0,y)=\psi(y)$ and $P(0,y)=\phi(y)$ admit a  projectively
flat Finsler metric with $\mathbf{K}=1$ near the origin if
and only if both $\psi$ and $\phi$ are real analytic on
$\R^n\setminus\{0\}$. The metric and its projective factor are unique
near the initial fibre.
\end{theoremA}

Our proof works directly with the $C^1$ first-order system and requires
neither homogeneity nor convexity for its local holomorphic
representation. For a smooth Finsler metric, it supplies the local
complex extensions used in Li's positive-curvature formula
\cite[Theorem~1.3]{Li}.

For $\mathbf{K}=-1$, the functions $Q_\pm=P\pm F$ satisfy real scalar
equations. For $\mathbf{K}=0$, the projective factor satisfies the same
equation. Li's implicit representations give local uniqueness for
smooth data in these two cases~\cite[Theorems~1.1 and~1.2]{Li},
while Theorem~\ref{thmA} gives local uniqueness in positive curvature.
We can therefore compare existing metrics on any connected coordinate
domain. The passage from local equality to equality on the whole
domain uses an open--closed argument and requires neither
completeness nor convexity.

\begin{theoremB}\label{thmB}
Let $F_1 = F_1(x,y)$ and $F_2 = F_2(x,y)$ be two smooth projectively flat
Finsler metrics on the same connected open set $\Omega\subset\R^n$, with projective factors $P_1 =P_1(x,y)$ and $P_2 =P_2(x,y)$ respectively, and
with the same constant flag curvature $\mathbf{K}$. Suppose that for some
$x_0\in\Omega$,
\[
F_1(x_0,y)=F_2(x_0,y),\qquad P_1(x_0,y)=P_2(x_0,y)
\quad\text{for all }y\in\R^n\setminus\{0\}.
\]
Then $F_1=F_2$ and $P_1=P_2$ on
$\Omega\times(\R^n\setminus\{0\})$.
\end{theoremB}

The explicit reconstruction from the initial pair can be found in
\cite[Theorems~1.1--1.3]{Li}; see
Theorem~\ref{thm:local-classification}. Thus the pair
\((\psi,\phi)=(F(x_0,\cdot),P(x_0,\cdot))\) determines \((F,P)\) near
the initial fibre.
Theorem~\ref{thmB} is a statement about solutions that already exist; it
does not claim that a metric extends to an arbitrary domain. It answers
the uniqueness questions for nonpositive curvature raised in
\cite[Remarks~3.2 and~4.2]{LZ}. Once local uniqueness on a full fibre is
available, the passage to a connected domain is an elementary open--closed
argument; the substance is the local uniqueness, which for $\mathbf K=1$
is Theorem~\ref{thmA}.

We then turn to the maximal connected domain of the real scalar
equation
\[
 Q_{x^k}=Q Q_{y^k}.
\]
The standard implicit formula $Q(x,y)=\alpha(y+xQ(x,y))$ gives a
local solution. We ask where it remains smooth in every nonzero
tangent direction.
Let $\alpha$ be a smooth real positively one-homogeneous
function on $\R^n\setminus\{0\}$. Its differential at $u\ne0$ is
\[
 \dd\alpha_u(v)
 =\left.\frac{\dd}{\dd t}\right|_{t=0}\alpha(u+tv)
 = \frac{\partial\alpha}{\partial y^j}(u)v^j.
\]
With $\Sn=\{u\in\R^n:|u|=1\}$, define
\begin{equation}\label{eq:D}
 \Delta_\alpha(x,\xi)=1-\dd\alpha_\xi(x),\qquad
 \calH_\alpha=\{x:\Delta_\alpha(x,u)>0\ \text{for every }u\in\Sn\}.
\end{equation}
Since $\dd\alpha$ is zero-homogeneous, the condition on unit vectors is
equivalent to the condition for all $\xi\ne0$.

\begin{theoremC}\label{thmC}
The set $\calH_\alpha$ is open and convex and contains the origin.
There is a unique smooth positively one-homogeneous solution $Q$ of
$Q_{x^k}=Q Q_{y^k}$ on
$\calH_\alpha\times(\R^n\setminus\{0\})$ with
$Q(0,\cdot)=\alpha$. It is determined by
\begin{equation}\label{eqQ_thmC}
 Q(x,y)=\alpha(y+xQ(x,y)).
\end{equation}
Moreover, if \(\Omega\subset\R^n\) is connected, open, contains the
origin, and a smooth positively one-homogeneous solution with
\(Q(0,\cdot)=\alpha\) exists on
\(\Omega\times(\R^n\setminus\{0\})\), then \(\Omega\subset\calH_\alpha\)
and the solution is the restriction of the one on \(\calH_\alpha\).
\end{theoremC}

Theorem~\ref{thmC} is proved in Subsection~\ref{subsec:maximal}. For each
$x\in\calH_\alpha$ the characteristic map
\[
 \Theta_x(\xi)=\xi-x\alpha(\xi)
\]
is a diffeomorphism of $\R^n\setminus\{0\}$ onto itself, with
\[
 \det\dd\Theta_x(\xi)=1-\dd\alpha_\xi(x)
\]
(see Corollary~\ref{cor:Tx}). The theorem determines the maximal
connected scalar domain containing $0$ without assuming positivity or
convexity of $\alpha$. Writing $\phi=P(0,\cdot)$ and $\psi=F(0,\cdot)$,
the corresponding domain is
\[
 \calH=
 \begin{cases}
 \calH_\phi,&\mathbf{K}=0,\\
 \calH_{\phi+\psi}\cap\calH_{\phi-\psi},& \mathbf{K}=-1.
 \end{cases}
\]

Over $\calH$ the reconstructed function $F$ is smooth and positive on
the slit bundle, but it defines a Finsler metric only on
\[
 \calO=\{x\in\calH:\tfrac12\Hess_y(F^2)(x,y)>0
                         \text{ for every }y\ne0\}.
\]
The component of $\calO$ containing $0$ is the maximal connected
regular extension of the initial metric in these coordinates.
The scalar solution can remain smooth beyond this component.
Angular Hessian tests describe $\calO$ in Section~\ref{sec:real}.
The straight-line invariance of $\Hess_y F$ gives planar convexity of
regularity components (Section~\ref{sec:boundary}). In dimensions
$n\ge3$ components need not be convex: the Hilbert-metric example of
Proposition~\ref{prop:nonconvex} has a connected nonconvex regularity
component. For equal Randers data the components are convex in every
dimension (Theorem~\ref{thm:randers}).

As an application we compute the strong-convexity locus for $\mathbf{K}=0$ and
equal Randers initial data
\[
 \phi(y)=\psi(y)=|y|+a\cdot y,\qquad |a|<1.
\]
Theorem~\ref{thm:randers} gives the locus in every dimension $n\ge2$,
together with the transition from one component to two convex
components; this extends the planar calculation of Li and
Zhao~\cite[Proposition~3.12]{LZ}. (Only the initial norm is assumed to
be of Randers type.)

In positive curvature our continuation criterion
(Section~\ref{sec:continuation}) assumes control of the holomorphic
branch, of its implicit denominator, and of the fundamental tensor;
analyticity in the interior is not enough to remove the
closed-hemisphere extension assumption in~\cite[Theorem~5.6]{LZ}.

The paper is organized as follows. Section~\ref{sec:prelim} fixes the
notation and recalls the local classification. Section~\ref{sec:holo}
proves the real analyticity of positive-curvature solutions and the exact
admissibility of their initial data (Theorem~\ref{thmA}), establishes
one-fibre uniqueness (Theorem~\ref{thmB}), and closes with a
variable-curvature example showing that the constant-curvature hypothesis
in Theorem~\ref{thmB} is essential. Section~\ref{sec:real} determines the
maximal scalar domain in nonpositive curvature (Theorem~\ref{thmC}),
reconstructs the metric on it, and reduces strong convexity to a test on
the initial data. Section~\ref{sec:continuation} gives the corresponding
conditional continuation criterion in positive curvature.
Section~\ref{sec:boundary} proves the straight-line invariance of the
fibre Hessian, deduces planar convexity and boundary tangency from it,
shows by example that components need not be convex for $n\ge3$
(Proposition~\ref{prop:nonconvex}), and computes the regularity locus
for equal Randers data in every dimension (Theorem~\ref{thm:randers}).

\section{Preliminaries}\label{sec:prelim}

This section fixes the notation and recalls the classical facts on which
the paper rests. For the general theory of Minkowski norms, Finsler
metrics, geodesics, and flag curvature we refer to
\cite{BaoChernShen2000}; only the formulas used below are recalled. We
begin with the conventions used throughout the paper.

Repeated indices are summed from $1$ to $n$, and $\delta^i_j$ is the
Kronecker symbol. The slit tangent bundle is
$T\Omega\setminus\{0\}=\Omega\times(\R^n\setminus\{0\})$.

For a function $f$ of a vector variable, $\dd f_u$ and $\Hess f_u$ are
its first and second differentials at $u$; for $F=F(x,y)$ we write
$\dd_yF$ and $\Hess_yF$ when $x$ is held fixed. A covector $\omega$ acts
by $\omega(v)=\omega_jv^j$, and $(v\otimes\omega)u=v\,\omega(u)$. For a
scalar $f$, $f_y=(f_{y^1},\dots,f_{y^n})$ is a row covector; for a
vector-valued map, $f_y$ is its Jacobian. We write
$\Dop=y^j\partial_{x^j}$.

Finally, for $y\ne0$, $\R y=\{ty:t\in\R\}$ is the radial line spanned by
$y$, and for $p\in\R^n$, $p+\R y=\{p+ty:t\in\R\}$ is the affine line
through $p$ in direction $y$. The quotient $\R^n/\R y$ has dimension
$n-1$.

\medskip

A scalar function $f$ is \emph{positively one-homogeneous} if
$f(ty)=tf(y)$ for $t>0$ and $y\ne0$; it is allowed to take either sign.

A \emph{Minkowski norm} on $\R^n$ is a continuous function
$\psi:\R^n\to[0,\infty)$, smooth on $\R^n\setminus\{0\}$, positively
one-homogeneous, with $\psi(y)=0$ only for $y=0$, and such that the matrix
$\bigl(\tfrac12(\psi^2)_{y^iy^j}\bigr)$ is positive definite at every
$y\ne0$. Reversibility is not assumed.

A \emph{Finsler metric} on $\Omega$ is a continuous function
$F=F(x,y):\Omega\times\R^n\to[0,\infty)$, smooth for $y\ne0$, whose
restriction $F(x,\cdot)$ to each tangent space is a Minkowski norm; its
\emph{fundamental tensor} is
\[
 g_{ij}=\tfrac12(F^2)_{y^iy^j},
 \qquad\text{with inverse } g^{ij}.
\]
Thus a Finsler metric is strongly convex by definition. We call a
solution \emph{regular} when its fundamental tensor is positive definite;
a scalar solution is not regular until this holds.

The geodesic coefficients $G^i$ of a Finsler metric $F$ are defined by
$\ddot x^i+2G^i(x,\dot x)=0$ for geodesics parametrized with constant
$F$-speed; locally
\[
 G^i=\frac14 g^{ij}\bigl([F^2]_{x^ky^j}y^k-[F^2]_{x^j}\bigr).
\]
A Finsler metric $F$ on a domain $\Omega\subset\R^n$ is
\emph{projectively flat} if the images of its geodesics are straight
lines. Equivalently, there are coordinates $(x^i)$ in which
\begin{equation}\label{eq:projective-G}
 G^i=Py^i
\end{equation}
for a positively one-homogeneous function $P$, the \emph{projective
factor} of $F$, given by
\begin{equation}\label{eq:projective-factor}
 P=\frac{\Dop F}{2F}.
\end{equation}
Hamel's characterization \cite{Hamel1903} states that $F$ is projectively
flat precisely when \eqref{eq:hamel} holds. Differentiating
$PF=\Dop F/2$ in $y^k$ and using \eqref{eq:hamel} gives
$F_{x^k}=(PF)_{y^k}$, the first equation of \eqref{eq:system}.

Following Berwald \cite{Berwald1929}, the curvature associated with the
coefficients $G^i$ is
\begin{equation}\label{eq:R}
 R^i{}_j=2G^i_{x^j}-y^kG^i_{x^ky^j}+2G^kG^i_{y^ky^j}
 -G^i_{y^k}G^k_{y^j},
\end{equation}
and for $G^i=Py^i$ it reduces to
\begin{equation}\label{eq:R-projective}
 R^i{}_j=(P^2-\Dop P)\delta^i_j
 +(2P_{x^j}-\Dop P_{y^j}-PP_{y^j})y^i .
\end{equation}
Constant flag curvature $\mathbf K$ means
$R^i{}_j=\mathbf K (F^2\delta^i_j-FF_{y^j}y^i)$. Comparing this with
\eqref{eq:R-projective} gives $\Dop P=P^2-\mathbf K F^2$, and
differentiating in $y^j$ yields the second equation of
\eqref{eq:system}. In particular
\begin{equation}\label{eq:flag-curvature}
 \mathbf K=\frac{P^2-\Dop P}{F^2}.
\end{equation}

The local classification of projectively flat Finsler metrics with constant flag curvature was solved by Shen \cite{Shen} (sufficiency) and Li \cite{Li} (necessity).

\begin{theorem}[\cite{Li,Shen}]
\label{thm:local-classification}
Let \(F = F(x,y)\) be a Finsler metric on an open neighbourhood \(\mathcal U\)
of the origin in \(\R^n\).
\begin{enumerate}
\item[(i)] \(F\) is projectively flat with \(\mathbf K=0\) if and only
if there are a Minkowski norm \(\psi\) and a positively one-homogeneous
function \(\phi\in C^\infty(\R^n\setminus\{0\})\) such that
\[
 F(x,y)=\psi(y+xP)\bigl(1+x^kP_{y^k}\bigr),
 \qquad P=\phi(y+xP).
\]

\item[(ii)] \(F\) is projectively flat with \(\mathbf K=-1\) if and
only if there are a Minkowski norm \(\psi\) and a positively
one-homogeneous function \(\phi\in C^\infty(\R^n\setminus\{0\})\) such
that
\[
 F(x,y)=\tfrac12\bigl(\Phi_+-\Phi_-\bigr),
 \qquad
 \Phi_\pm=(\phi\pm\psi)(y+x\Phi_\pm).
\]

\item[(iii)] Suppose that \(F\) extends to \(F(x,y+xz)\) for
\(z\in\C \). Then \(F\) is projectively flat with \(\mathbf K=1\) if
and only if there are a Minkowski norm \(\psi\) and a positively
one-homogeneous function \(\phi\in C^\infty(\R^n\setminus\{0\})\) such
that
\[
 F(x,y)=\im \Psi(x,y),
 \qquad
 \Psi=(\phi+i\psi)(y+x\Psi).
\]
\end{enumerate}
In each case, \(\psi(y)=F(0,y)\) and \(\phi(y)=P(0,y)\).
\end{theorem}

Thus the local problem is determined by the pair $(\psi,\phi)$ on one
fibre. In (iii) the parameter $z$ is a complex scalar and $y+xz\in\C^n$;
the representation holds for $y$ near a real slit direction and $z$ near
$Q(x,y)$. For $\mathbf K=1$ condition (iii) is satisfied automatically:
by Theorem~\ref{thmA}, every smooth projectively flat metric of flag
curvature one is real analytic, and hence extends holomorphically to a
neighbourhood of each point $y+xQ(x,y)$ with $y\ne0$.

The construction of solutions for $\mathbf{K} \le0$ hinges on the existence of a
unique root of the implicit equation $q=\alpha(y+xq)$. The following is
the local form of this fact, proved in \cite[Lemma~3.2 and Theorem~3.3]{Li}.

\begin{lemma}[\cite{Li}]\label{lem:unique-root}
Let $\alpha\in C^\infty(\R^n\setminus\{0\})$ be positively
one-homogeneous, extended by $\alpha(0)=0$. There is $\varepsilon>0$ such
that for every $x$ with $|x|<\varepsilon$ and every $y\ne0$, the function
$h(q)=q-\alpha(y+xq)$ has exactly one zero $q=Q(x,y)$; moreover $Q$ is
$C^\infty$ near $(x,y)$ and  positively one-homogeneous in $y$.
\end{lemma}

The global version, with the ball $|x|<\varepsilon$ replaced by the
condition $x\in\calH_\alpha$, is part of
Theorem~\ref{thmC}.

While the local existence for nonpositive curvature relies on the uniqueness of the real implicit root (Lemma~\ref{lem:unique-root}), the positive curvature case (\(\mathbf{K} =1\)) hinges on the analogous unique holomorphic root. This fundamental existence and uniqueness result for the complex implicit equation was established by Li \cite[Theorem 6.1]{Li}.
That construction required the initial data to be complexifiable (i.e., \emph{a priori} real analytic). However, thanks to the automatic analyticity established in Theorem~\ref{thmA} (Corollary~\ref{cor:analytic}), this condition is naturally satisfied for any smooth positive-curvature metric. We recall the lemma here in our current notation.

\begin{lemma}[\cite{Li}]\label{lem:holo-root}
Let \(\alpha\) be a real analytic, positively one-homogeneous function on \(\R^n\setminus\{0\}\), and let \(h\) be its holomorphic extension near a point \(u\in\Sn\). Consider the complex implicit equation
\[
 q = h(y+xq).
\]
Near \((x,y,q)=(0,u,\alpha(u))\), this equation has a unique holomorphic solution \(q = Q(x,y)\).
\end{lemma}

\medskip

We record three facts used repeatedly below. The last one computes the
Hessian of an implicit branch $Q=\alpha(y+xQ)$; it is applied to
$\alpha_\pm=\phi\pm\psi$ in Subsection~\ref{subsec:convex} and to the holomorphic
branch in Section~\ref{sec:continuation}.

\begin{lemma}\label{lem:rank-one}
Let $x\in\R^n$ (or $\C^n$) and let $\omega$ be a covector with
$\omega(x)\ne1$, and put $A=\Id-x\otimes\omega$. Then
\[
 A^{-1}=\Id+\frac{x\otimes\omega}{1-\omega(x)},\qquad
 \det A=1-\omega(x),\qquad
 A^{-1}x=\frac{x}{1-\omega(x)}.
\]
\end{lemma}

\begin{proof}
A direct multiplication gives
\[
 (\Id-x\otimes\omega)\Bigl(\Id+\frac{x\otimes\omega}{1-\omega(x)}\Bigr)
 =\Id+\frac{x\otimes\omega}{1-\omega(x)}-x\otimes\omega
 -\frac{\omega(x)}{1-\omega(x)}x\otimes\omega
 =\Id,
\]
using $(x\otimes\omega)(x\otimes\omega)=\omega(x)\,x\otimes\omega$.
Hence the displayed formula for $A^{-1}$ holds. On the hyperplane
$\ker\omega$ the map $A$ is the identity, and on the line $\R x$ it acts
by $1-\omega(x)$; therefore $\det A=1-\omega(x)$. Finally,
\[
 A^{-1}x=\Bigl(\Id+\frac{x\otimes\omega}{1-\omega(x)}\Bigr)x
 =x+\frac{\omega(x)}{1-\omega(x)}x=\frac{x}{1-\omega(x)}.
\]
\end{proof}

The following criterion is standard; see \cite{BaoChernShen2000}. We
recall the short proof because the quotient formulation is used
repeatedly below.

\begin{lemma} [\cite{BaoChernShen2000}]
\label{lem:angular}
Let \(f\) be smooth, positive, and positively one-homogeneous on
\(\R^n\setminus\{0\}\). For \(y\ne0\), write
\[
g_y:=\tfrac12\Hess(f^2)_y
=\dd f_y\otimes\dd f_y+f(y)\Hess f_y
\]
for the fundamental tensor of \(f\) at \(y\). Then
\(\Hess f_y(y,\cdot)=0\). Hence \(\Hess f_y\) induces a symmetric
bilinear form \(\bar H_y\) on the quotient \(\R^n/\R y\), defined by
\[
\bar H_y([v],[w]):=\Hess f_y(v,w),
\]
independently of the chosen representatives. We call \(\bar H_y\)
the \emph{angular Hessian} of \(f\) at \(y\). The fundamental tensor
\(g_y\) is positive definite on \(\R^n\) if and only if \(\bar H_y\)
is positive definite on \(\R^n/\R y\).
\end{lemma}

\begin{proof}
Since \(\dd f_y(y)=f(y)>0\), the hyperplane
\(E_y=\ker \dd f_y\) is a complement of \(\R y\). From the displayed
identity,
\[
g_y(y,y)=f(y)^2>0,\qquad
g_y(y,v)=f(y)\dd f_y(v)=0\quad(v\in E_y),
\]
and
\[
g_y(v,w)=f(y)\Hess f_y(v,w)\quad(v,w\in E_y).
\]
Thus \(\R^n=\R y\oplus E_y\) is an orthogonal decomposition for
\(g_y\), and \(g_y\) is positive on \(\R y\). Therefore \(g_y\) is
positive definite on \(\R^n\) if and only if \(\Hess f_y\) is
positive definite on \(E_y\). The quotient map restricts to an
isomorphism \(E_y\to\R^n/\R y\), so this is equivalent to positive
definiteness of \(\bar H_y\).
\end{proof}

\begin{lemma}\label{lem:impHess}
Let $Q$ be a $C^1$ function, positively one-homogeneous in $y$, on an
open subset of $\R^n_x\times(\R^n_y\setminus\{0\})$, satisfying
\[
 Q(x,y)=\alpha\bigl(y+xQ(x,y)\bigr),\qquad y\ne0 .
\]
Assume either
\begin{enumerate}
\item[(i)] $\alpha\in C^\infty(\R^n\setminus\{0\})$ is real-valued,
positively one-homogeneous, extended by $\alpha(0)=0$, and $Q$ is
real-valued; or
\item[(ii)] $\alpha=h$ is holomorphic and positively one-homogeneous on
an open set $W\subset\C^n$, $Q$ is complex-valued, and
$y+xQ(x,y)\in W$ for every $(x,y)$ in the domain of $Q$.
\end{enumerate}
Put
\[
 \xi=y+xQ(x,y),\qquad \Delta=1-\dd\alpha_\xi(x),\qquad
 A=\Id-x\otimes\dd\alpha_\xi .
\]
Then $\xi\ne0$, $\Delta\ne0$, $Q$ is in fact $C^\infty$, and
\begin{equation}\label{eq:impHess}
 \Hess_yQ=\frac{A^{-T}(\Hess\alpha)_\xi A^{-1}}{\Delta},
\end{equation}
where the Hessian is real bilinear in case (i) and complex bilinear in
case (ii).
\end{lemma}

\begin{proof}
The computation below is formal and applies to both cases. In (i) all
differentials are real-linear and the Hessians real bilinear; in (ii) they
are complex-linear and complex bilinear.

If $\xi=0$, then in case (i) $\alpha(0)=0$ by definition, and in case
(ii) $0=\xi\in W$ and $h(0)=0$ by positive homogeneity. Hence
$Q=\alpha(0)=0$, and then $y=-xQ=0$, contrary to $y\ne0$.

Differentiating $Q=\alpha(\xi)$ and $\xi=y+xQ$ gives
\[
 \dd Q=\dd\alpha_\xi\circ\dd\xi,\qquad \dd\xi=\Id+x\otimes\dd Q .
\]
Substituting the second identity into the first,
\[
 \dd Q=\dd\alpha_\xi+\dd\alpha_\xi(x)\,\dd Q,
 \qquad\text{so}\qquad \Delta\,\dd Q=\dd\alpha_\xi .
\]
In particular $\Delta\ne0$: otherwise $\dd\alpha_\xi=0$, so
$\dd\alpha_\xi(x)=0$ and $\Delta=1$. Hence
\[
 \dd\xi=\Id+\frac{x\otimes\dd\alpha_\xi}{\Delta}=A^{-1},
\]
and  the implicit function theorem, applied with $\Delta\ne0$, makes
$Q$ of class $C^\infty$.

Differentiating $\dd Q=\dd\alpha_\xi\circ\dd\xi$ once more,
\[
 \Hess Q(u,v)=(\Hess\alpha)_\xi(A^{-1}u,A^{-1}v)
 +\dd\alpha_\xi\bigl(\dd^2\xi(u,v)\bigr).
\]
Since $\dd^2\xi(u,v)=x\,\Hess Q(u,v)$ and $\dd\alpha_\xi(x)=1-\Delta$,
\[
 \Hess Q(u,v)=(\Hess\alpha)_\xi(A^{-1}u,A^{-1}v)
 +(1-\Delta)\Hess Q(u,v),
\]
hence $\Delta\,\Hess Q(u,v)=(\Hess\alpha)_\xi(A^{-1}u,A^{-1}v)$, which is
\eqref{eq:impHess} as bilinear forms. This proves \eqref{eq:impHess}; the Hessian is real bilinear in case (i)
and complex bilinear in case (ii).
\end{proof}

The lemma is conditional: it computes the Hessian of a solution but does
not produce one. In the real case existence is supplied by
Theorem~\ref{thmC}, which for smooth positively one-homogeneous
$\alpha$ constructs the solution of $Q=\alpha(y+xQ)$ on
$\calH_\alpha$ (its proof incorporates the unique-root argument of
\cite[Lemma~3.2 and Theorem~3.3]{Li}). This is the form used in
Subsection~\ref{subsec:convex}, with $\alpha_\pm=\phi\pm\psi$, where the Hessians
of $Q_\pm$ give
$\Hess_yF=\frac12(\Hess_yQ_+-\Hess_yQ_-)$. In the holomorphic case
existence near each point is supplied by \cite[Theorem~6.1]{Li}, recalled
as Lemma~\ref{lem:holo-root} in Section~\ref{sec:prelim}; along a path it is part
of the hypothesis of Theorem~\ref{thm:continue}. This is the form used in
Section~\ref{sec:continuation}: the formula exhibits $\Hess_yQ$ as a continuous
function of $(x,u,Q)$, which is what lets strong convexity pass to the
limit.

\section{Analyticity, initial data, and one-fibre uniqueness}\label{sec:holo}
In the positive curvature case (\(\mathbf{K}=1\)), the implicit constructions of Shen~\cite{Shen} and Li~\cite{Li} rely on a holomorphic representation of the solution. This section establishes the two main results for this case: Theorem~\ref{thmA} (automatic real analyticity and the exact admissibility of initial data) and Theorem~\ref{thmB} (global one-fibre uniqueness on an arbitrary connected domain).

First, we prove that smoothness alone forces the metric to be real analytic (Corollary~\ref{cor:analytic}). Second, leveraging this analyticity and the unique holomorphic root from Lemma~\ref{lem:holo-root}, we establish the exact smooth admissibility condition for the initial data (Theorem~\ref{thm:data}). Together, these two steps yield a complete proof of Theorem~\ref{thmA}. Finally, the local uniqueness is propagated to any connected coordinate domain via a topological open-closed argument, yielding Theorem~\ref{thmB}.
We close with a variable-curvature example (Proposition~\ref{prop:counter}) showing that the
constant-curvature hypothesis in Theorem~\ref{thmB} cannot be dropped.

Throughout Subsections~\ref{subsec:Local holomorphic representation}--\ref{subsec:Initial data and one-fibre uniqueness}, let \(F=F(x,y)\) be a smooth projectively flat Finsler metric of constant flag curvature \(\mathbf{K}=1\), and let \(P=P(x,y)\) be its projective factor. The restriction of these functions to the initial tangent fibre at the origin defines the positively one-homogeneous initial data
\[
 \phi(y)=P(0,y), \qquad \psi(y)=F(0,y).
\]

\subsection{Local holomorphic representation}\label{subsec:Local holomorphic representation}
The following argument uses only a first-order equation and does not require homogeneity.

\begin{theorem}\label{thm:holo}
Let \(Q\) be a complex-valued \(C^1\) function on an open subset of \(\R_x^n\times\R_y^n\), satisfying
\begin{equation}\label{eq:burgers}
 Q_{x^k}=Q Q_{y^k},\qquad k=1,\ldots,n.
\end{equation}
Suppose \(\im Q(x_0,y_0)\ne0\). There is a holomorphic function \(h\) near \(y_0\in\C^n\) such that, near \((x_0,y_0)\),
\begin{equation}\label{eq:local-implicit}
 Q(x,y)=h\big(y+(x-x_0)Q(x,y)\big).
\end{equation}
In particular, \(Q\) is real analytic near \((x_0,y_0)\). On the real slice, \(h(y)=Q(x_0,y)\).
\end{theorem}

\begin{proof}
We first construct a real change of coordinates. Regard $x$ and $y$ as
column vectors, put $X=x-x_0$ and $q_0=Q(x_0,y_0)$, and define
\[
 \Phi(x,y):=y+XQ(x,y)\in\C^n .
\]
Since $\dd X=\dd x$, differentiating $\Phi=y+XQ$ gives
\[
 \dd\Phi=\dd y+Q\,\dd x+X\,\dd Q .
\]
Let $(\delta x,\delta y)\in\R^n\times\R^n$ be an arbitrary increment at
$(x_0,y_0)$. At that point $X=0$, so the term $X\,\dd Q$ drops out and
\[
 \dd\Phi_{(x_0,y_0)}(\delta x,\delta y)=\delta y+q_0\,\delta x .
\]
If this vector vanishes, its real and imaginary parts vanish separately:
\[
 \delta y+(\re q_0)\delta x=0,
 \qquad
 (\im q_0)\delta x=0 .
\]
Since $\im q_0=\im Q(x_0,y_0)\ne0$ by hypothesis, the second equation
gives $\delta x=0$, and then the first gives $\delta y=0$. Thus
$\dd\Phi_{(x_0,y_0)}$ is an injective real-linear map
$\R^n_x\times\R^n_y\to\C^n$. Both spaces have real dimension $2n$, so the
map is an isomorphism and the real Jacobian of $\Phi$ at $(x_0,y_0)$ is
invertible. Since $Q$ is $C^1$, so is $\Phi$; the real inverse function
theorem therefore gives neighbourhoods $U$ of $(x_0,y_0)$ and $V$ of
$y_0$ such that $\Phi:U\to V$ is a $C^1$ diffeomorphism. We write
$z=\Phi(x,y)$.

We next show that $Q$, expressed in the coordinates $z$, is holomorphic.
Let $r=(Q_{y^1},\ldots,Q_{y^n})$ be the complex row vector of fibre
derivatives, and let $\theta=\dd y+Q\,\dd x$ be a column of
complex-valued differential forms. By \eqref{eq:burgers},
\[
 \dd Q
 =Q_{x^k}\dd x^k+Q_{y^k}\dd y^k
 =Q_{y^k}\bigl(\dd y^k+Q\,\dd x^k\bigr)
 =r\,\theta .
\]
Consequently
\[
 \dd z=\dd y+Q\,\dd x+X\,\dd Q
 =\theta+X\,\dd Q
 =(\Id+X\otimes r)\,\theta
 =B\,\theta ,
 \qquad B=\Id+X\otimes r ,
\]
where $(X\otimes r)^j{}_k=X^j r_k$. Since $B=\Id$ at the centre, we may
shrink $U$ so that $B$ is invertible on $U$; by
Lemma~\ref{lem:rank-one},
\[
 B^{-1}=\Id-\frac{X\otimes r}{1+r(X)} .
\]
Hence $\theta=B^{-1}\dd z$ and therefore
\begin{equation}\label{eq:dQ-dz}
 \dd Q=rB^{-1}\dd z .
\end{equation}

Define $h=Q\circ\Phi^{-1}$ on $V$, so that $h(z)=Q(x,y)$ whenever
$z=\Phi(x,y)$. By the chain rule for differentials,
\[
 \dd h_z=\dd Q_{(x,y)}\circ\bigl(\dd\Phi_{(x,y)}\bigr)^{-1},
 \qquad (x,y)=\Phi^{-1}(z).
\]
Substituting \eqref{eq:dQ-dz} and using $\dd z=\dd\Phi_{(x,y)}$,
\[
 \dd h_z
 =(rB^{-1})(x,y)\circ\dd z_{(x,y)}\circ
 \bigl(\dd\Phi_{(x,y)}\bigr)^{-1}
 =(rB^{-1})(x,y).
\]
Thus, in the $z$-coordinates, $\dd h$ is multiplication by the complex row
vector $(rB^{-1})(\Phi^{-1}(z))$:
\[
 \dd h=\bigl[(rB^{-1})\circ\Phi^{-1}\bigr](z)\,\dd z .
\]
Writing $\dd z=(\dd z^1,\ldots,\dd z^n)$ and
$\dd\bar z=(\dd\bar z^1,\ldots,\dd\bar z^n)$, this says that $\dd h$ has
no $\dd\bar z$ component. Hence all Wirtinger derivatives
$h_{\bar z^j}$ vanish; since $h$ is $C^1$, it is holomorphic.

Finally, we recover analyticity in the original variables. By
construction $h(y+XQ)=Q$, and $h(y)=Q(x_0,y)$ when $X=0$. For complex
variables $(X,y,q)$ near $(0,y_0,q_0)$, set $E(X,y,q)=q-h(y+Xq)$. This
function is holomorphic and vanishes at $(0,y_0,q_0)$; moreover
\[
 E_q(X,y,q)=1-\sum_{j=1}^n X^j h_{z^j}(y+Xq),
 \qquad E_q(0,y_0,q_0)=1 .
\]
The holomorphic implicit function theorem gives a unique holomorphic
solution $q=q(X,y)$ near $(0,y_0)$. The original $C^1$ function
$Q(x_0+X,y)$ satisfies $Q=h(y+XQ)$ and, by continuity, is close to $q_0$
after $U$ is reduced. Implicit uniqueness identifies it with the
restriction of $q(X,y)$ to real $(X,y)$. Hence $Q$ is real analytic near
$(x_0,y_0)$.
\end{proof}

\begin{corollary}\label{cor:analytic}
Every smooth projectively flat Finsler metric of constant flag curvature $\mathbf{K}=1$ is real analytic in \((x,y)\) in its projective coordinates, and so is its projective factor.
\end{corollary}

\begin{proof}
With \(Q=P+iF\), equations~\eqref{eq:system} give
\[
 \begin{aligned}
 Q_{x^k} &=PP_{y^k}-FF_{y^k}+i(PF_{y^k}+FP_{y^k})\\ &=(P+iF)(P_{y^k}+iF_{y^k})=Q Q_{y^k}.
 \end{aligned}
\]
Theorem~\ref{thm:holo} applies at every point \((x,y)\) with \(y\ne0\). Taking real and imaginary parts proves that \(P\) and \(F\) are real analytic there.
\end{proof}

\subsection{Initial data and one-fibre uniqueness}\label{subsec:Initial data and one-fibre uniqueness}
Local uniqueness for \(\mathbf K=1\) follows from
Corollary~\ref{cor:analytic} and Li's representation
\cite[Theorem~1.3]{Li}. For the proof of Theorem~\ref{thm:data}, we use
Lemma~\ref{lem:holo-root} to construct the holomorphic branch; the same
computation yields the derivative formulas needed in
Section~\ref{sec:continuation}.

\begin{theorem}\label{thm:data}
Let \(\psi\) be a smooth Minkowski norm, and let \(\phi\) be smooth and positively one-homogeneous on \(\R^n\setminus\{0\}\). The initial data \(F(0,y)=\psi(y)\) and \(P(0,y)=\phi(y)\) admit a smooth projectively flat Finsler metric of flag curvature \(\mathbf{K}=1\) near \(0\) if and only if both functions are real analytic on the slit space. The metric and its projective factor are unique as germs over the entire initial fibre.
\end{theorem}
\begin{proof}
Necessity follows from Corollary~\ref{cor:analytic}.

For sufficiency, fix \(u\in\Sn\) and complexify \(\alpha=\phi+i\psi\) near \(u\). Denote the resulting holomorphic germ by \(h\). By Lemma~\ref{lem:holo-root}, the equation
\begin{equation}\label{eq:data-implicit}
 q=h(y+xq)
\end{equation}
has a unique holomorphic solution \(Q(x,y)\) near \((0,u)\).
The holomorphic extensions agree near each real point of overlap:
their difference is holomorphic there and vanishes on an open subset
of the real slice, so all its real derivatives, and hence all its
holomorphic derivatives, vanish; it is therefore identically zero near
that point. After shrinking, implicit uniqueness makes the local
solutions agree on overlaps.

Covering \(\Sn\) by finitely many such neighbourhoods, we get a common ball \(|x|<\varepsilon\) on which these solutions define \(Q\) for real \(y\) in an annular neighbourhood of the sphere. To extend \(Q\) homogeneously to all \(y\ne0\), we check compatibility with these annular solutions. Euler's identity \(y^j\alpha_{y^j}=\alpha\) complexifies to
\begin{equation}\label{eq:euler-h}
 z^j h_{z^j}=h.
\end{equation}
Equivalently, the germ of \(h\) scales by degree one under small positive real scalings. Uniqueness in~\eqref{eq:data-implicit} then gives \(Q(x,ty)=tQ(x,y)\) whenever both sides are in the annulus, and hence the radial extensions agree.

To extract the derivative formulas required in Section~\ref{sec:continuation}, we differentiate the implicit equation. Let \(\omega =\dd h_z = (\omega_1, \cdots, \omega_n)\), where \( \omega_k = \partial h /\partial z_k(z) \), \(z=y+xQ\), and put \(\Delta=1-\omega(x)\). Then
\begin{equation}\label{eq:firstderivatives}
 Q_{y^k}=\frac{\omega_k}{\Delta},\qquad Q_{x^k}=\frac{Q \omega_k}{\Delta}=Q Q_{y^k}.
\end{equation}
Set \(P=\re Q\) and \(F=\im Q\). The real and imaginary parts of \eqref{eq:burgers} give the system \eqref{eq:system}. At \(x=0\), \(F=\psi>0\) and \(\frac12\Hess_y(F^2)=\frac12\Hess(\psi^2)>0\). By continuity and compactness of the unit sphere, positivity and positive definiteness persist on a smaller ball, and homogeneity extends them to all nonzero directions.

Finally, to prove uniqueness, suppose two smooth solutions have the same initial data on the entire fibre. By Theorem~\ref{thm:holo}, each is represented by a holomorphic germ whose restriction to the real slice is the given initial data \(\alpha\). Since the difference vanishes on the real slice, all its real derivatives vanish there. The holomorphic derivatives are linear combinations of these, so all Taylor coefficients at each point of the slice are zero. Hence the difference is identically zero near the slice, and the two germs coincide. Consequently, both solutions satisfy the same implicit equation \(q=h(y+xq)\), and the uniqueness of the holomorphic root is guaranteed by Lemma~\ref{lem:holo-root}. Hence \(Q\), and therefore \((P,F)\), is uniquely determined near the initial fibre.
\end{proof}

\begin{proof}[Proof of Theorem~\ref{thmA}]
The necessity of real analyticity is established by Corollary~\ref{cor:analytic}. The sufficiency of the initial data and the local uniqueness of the solution are established by Theorem~\ref{thm:data}. Combining these two results completes the proof.
\end{proof}

\begin{proposition}\label{prop:nonanalytic}
There are smooth strongly convex initial norms, arbitrarily $C^2$-close to the Euclidean norm on the unit sphere, for which no smooth positive-constant-curvature solution with those initial data exists, even locally in the base.
\end{proposition}

\begin{proof}
Define $b(s)=0$ for $s\le0$ and $b(s)=e^{-1/s^2}$ for $s>0$, and set
\[
 \psi_\epsilon(y)=|y|\left[1+\epsilon\, b\left(\frac{y^1}{|y|}\right)\right],
 \qquad \phi(y)=0.
\]
The functions $\psi_\epsilon$ are smooth and positively homogeneous, and
they converge to $|y|$ in $C^2$ on the unit sphere as $\epsilon\to0$; so
for small $\epsilon>0$ positivity and the positive definiteness of
$\frac12\Hess(\psi_\epsilon^2)$ persist uniformly.

Let $p$ be a nonzero point with $p^1=0$.
Since $b$ is flat at $0$ and $b(0)=0$, and since
$y\mapsto y^1/|y|$ is smooth near $p$, the perturbation $\epsilon\,|y|\,b\left(\frac{y^1}{|y|}\right)$ vanishes to
infinite order at $p$: all its derivatives at $p$ are zero. Hence the
Taylor series of $\psi_\epsilon$ at $p$ is that of $|y|$, which converges
to $|y|$ in a neighbourhood of $p$.
If $\psi_\epsilon$ were real analytic at $p$, it would equal $|y|$ there.
But $b\bigl(q^1/|q|\bigr)>0$ whenever $q^1>0$, and therefore
$\psi_\epsilon(q)>|q|$ at such $q$, a contradiction.
Thus $\psi_\epsilon$ is not real analytic at $p$,
and Theorem~\ref{thm:data} excludes a solution of constant positive flag
curvature with these initial data.
\end{proof}

\begin{proof}[Proof of Theorem~\ref{thmB}]
After the common rescaling $F_i\mapsto\sqrt{|\mathbf{K}|}F_i$ when $\mathbf{K}\ne0$,
we may assume $\mathbf{K}\in\{-1,0,1\}$. This leaves the projective factors
and the initial agreement unchanged. Translate $x_0$ to the origin.

The proof uses the following local fact, applied after translating
the comparison point to the origin. If two solutions of
\eqref{eq:system} with the same $\mathbf{K}$ agree on the fibre over a point,
then they agree on $U\times(\R^n\setminus\{0\})$ for some neighbourhood
$U$ of that point, and $U$ can be chosen independently of the fibre
direction. For $\mathbf{K}=1$, this is the uniqueness part of
Theorem~\ref{thm:data}.
For $\mathbf{K}=0$, both projective factors satisfy
$P=\phi(y+xP)$ with the same $\phi$; by Lemma~\ref{lem:unique-root}, the
root is unique near its initial value, so the factors coincide,
 and the reconstruction formula
\[
F=\psi(y+xP)\bigl(1+x^kP_{y^k}\bigr)
\]
gives the same metric function.

For $\mathbf{K}=-1$, the functions $Q_\pm=P\pm F$ satisfy
$Q_\pm=(\phi\pm\psi)(y+xQ_\pm)$. By the local classification
(Theorem~\ref{thm:local-classification}; see
\cite[Theorem~1.2]{Li}), the branch with prescribed initial data
is unique, so the corresponding branches of the two metrics coincide;
hence $P=\tfrac12(Q_++Q_-)$ and $F=\tfrac12(Q_+-Q_-)$ coincide as well.

Let
\[
S=\{x\in\Omega:F_1(x,y)=F_2(x,y)\ \text{and}\ P_1(x,y)=P_2(x,y)
\quad\text{for all }y\ne0\}.
\]
Then $x_0\in S$. By the local fact, each point of $S$ has a
neighbourhood contained in $S$, so $S$ is open in $\Omega$. If
$x_k\in S$ and $x_k\to x_*\in\Omega$, then for every fixed $y\ne0$
\[
F_1(x_*,y)=\lim_{k\to\infty}F_1(x_k,y)
=\lim_{k\to\infty}F_2(x_k,y)=F_2(x_*,y),
\]
and in the same way $P_1(x_*,y)=P_2(x_*,y)$; thus $x_*\in S$ and $S$ is
relatively closed. As $\Omega$ is connected and $S$ is nonempty,
$S=\Omega$.
\end{proof}

\subsection{Sharpness of the constant-curvature hypothesis}
\label{subsec:sharpness}

The hypothesis of constant curvature in Theorem~\ref{thmB} is essential. The
following family shows this in every dimension \(n\ge2\).

\begin{proposition}\label{prop:counter}
For $\varepsilon\geq0$, define on $T\R^n\setminus\{0\}$
\begin{equation}\label{eq:counter}
 F_\varepsilon := F_\varepsilon(x,y)=|y|+\varepsilon\left(|x|^2|y|+\frac{\langle x,y\rangle^2}{|y|}\right).
\end{equation}
These $F_\varepsilon(x,y)$ are real analytic, reversible, projectively flat, complete Finsler metrics.
For every $\varepsilon$,
\[
 F_\varepsilon(0,y)=|y|,\qquad P_\varepsilon(0,y)=0.
\]
For $\varepsilon>0$ the flag curvature is not constant. Thus one-fibre initial data do not determine a general projectively flat Finsler metric, even under analyticity, reversibility, and completeness.
\end{proposition}

\begin{proof}
Positivity, reversibility, homogeneity and analyticity on the slit bundle follow from the formula. Put $r=|y|$, $s=\langle x,y\rangle$. For $w\perp y$, differentiation gives
\[
 (\Hess_yF_\varepsilon)(w,w)=\frac1r\left[
 \left(1+\varepsilon\left(|x|^2-\frac{\langle x,y\rangle^2}{|y|^2}\right)\right)|w|^2
 +2\varepsilon\langle x,w\rangle^2\right]>0
\]
for $w\ne0$, and strong convexity follows from Lemma~\ref{lem:angular}.

Hamel's equation holds:
\[
 (F_{\varepsilon})_{x^k}= 2 \varepsilon \left(  x^k |y|+ \frac{\langle x,y\rangle y^k}{|y|} \right) =y^j (F_{\varepsilon})_{x^j y^k}.
\]
Thus $F_\varepsilon$ is projectively flat. Its projective factor is
\begin{equation}\label{eq:P-counter}
 P_\varepsilon(x,y)=\frac{2\varepsilon\langle x,y\rangle|y|}{F_\varepsilon(x,y)}.
\end{equation}
The two identities for the initial data follow at once.

Since \(F_\varepsilon\ge|y|\), an \(F_\varepsilon\)-unit-speed geodesic
has Euclidean speed at most one. Hence, on any finite time interval,
its base points remain in a compact set \(\mathcal{C}\). By homogeneity and the
compactness of \(\mathcal{C}\times\Sn\), there are constants \(0<c_1 < c_2\) such
that
\[
 c_1 |y|\le F_\varepsilon(x,y)\le c_2|y|
 \qquad(x\in \mathcal{C},\ y\ne0).
\]
Since \(F_\varepsilon(\gamma,\dot\gamma)=1\), the velocity lies in the
compact annulus
\[
 \frac{1}{c_2}\le|\dot\gamma|\le\frac{1}{c_1}.
\]
The geodesic equation is smooth on \(\mathcal{C}\times\{y:1/c_2 \le|y|\le1/c_1\}\),
so the solution extends beyond every finite time. Reversibility gives
completeness in both directions.

Finally, we show that the flag curvature is not constant. Put
\(r=|y|\) and \(s=\langle x,y\rangle\). Then
\[
F_\varepsilon
=r+\varepsilon\left(|x|^2r+\frac{s^2}{r}\right),
\qquad
P_\varepsilon=\frac{2\varepsilon sr}{F_\varepsilon}.
\]
Since \(\Dop s=r^2\) and \(\Dop F_\varepsilon=4\varepsilon sr\),
a direct computation gives
\[
\Dop P_\varepsilon
=\frac{2\varepsilon r^2\left(rF_\varepsilon-4\varepsilon s^2\right)}
{F_\varepsilon^2}.
\]
Substituting this and
\(P_\varepsilon^2=4\varepsilon^2s^2r^2/F_\varepsilon^2\) into
\(\mathbf K=(P_\varepsilon^2-\Dop P_\varepsilon)/F_\varepsilon^2\),
we obtain
\[
\mathbf K(x,y)
=\frac{2\varepsilon |y|^2\left(
5\varepsilon\langle x,y\rangle^2
-|y|^2(1+\varepsilon|x|^2)\right)}
{F_\varepsilon(x,y)^4}.
\]
Hence
\(\mathbf K\) is not constant for \(\varepsilon>0\).
\end{proof}

\section{Real characteristics and nonpositive curvature}\label{sec:real}
For $\mathbf K\le0$ the complex change of variables used for $\mathbf K=1$
is not available, and the scalar reductions of \eqref{eq:system} satisfy
the real first-order equation
\begin{equation}\label{eq:scalar}
 Q_{x^k}=Q\,Q_{y^k},\qquad 1\le k\le n .
\end{equation}
Fix a smooth positively one-homogeneous function $\alpha$ on
$\R^n\setminus\{0\}$, extended by $\alpha(0)=0$. The problem of this
section is to find the largest connected open set $\Omega\subset\R^n$
containing $0$ on which \eqref{eq:scalar} has a $C^\infty$ solution $Q$,
positively one-homogeneous in $y$, with $Q(0,\cdot)=\alpha$. The answer is
the set $\calH_\alpha$ of \eqref{eq:D}; this is Theorem~\ref{thmC}.

The argument has three parts, carried out in
Subsections \ref{subsec:chars}--\ref{subsec:maximal}. First, \eqref{eq:scalar} is
equivalent to the constancy of $Q$ along the characteristic curves of the
equation, and these curves are complete over every finite base path: the
one starting from $\xi$ over the origin ends at
$\Theta_x(\xi)=\xi-x\alpha(\xi)$. Second, for a solution on $\Omega$ the
endpoint map $\Theta_x$ is a diffeomorphism of the slit fibre for every
$x\in\Omega$; since $\det\dd\Theta_x(\xi)=\Delta_\alpha(x,\xi)$, this
makes $\Delta_\alpha>0$ necessary, so $\Omega\subset\calH_\alpha$. Third,
on $\calH_\alpha$ the implicit equation $Q=\alpha(y+xQ)$ has a unique
smooth solution, so the condition is also sufficient. At the end of
Subsection~\ref{subsec:maximal}, Proposition~\ref{prop:line} and
Corollary~\ref{cor:line} treat base domains containing a complete affine
line: the initial data are then linear, and for constant negative
curvature no such domain exists.

Once $Q$ is known, we recover the projective factor $P$ and the metric
function $F$ in Subsection~\ref{subsec:reconstruction}: for $\mathbf K=0$ from $P$ alone,
and for $\mathbf K=-1$ from $Q_\pm=P\pm F$. Finally, Subsection~\ref{subsec:convex} determines
the subset $\calO\subset\calH_\alpha$ on which $F$ is strongly convex, and
hence a Finsler metric.

\subsection{Characteristic curves}\label{subsec:chars}
Throughout this subsection $Q\in C^\infty(\Omega\times(\R^n\setminus\{0\}))$
is positively one-homogeneous in $y$; it need not satisfy \eqref{eq:scalar}.
For a piecewise $C^1$ path $\gamma:[0,1]\to\Omega$ and $\xi\in\R^n\setminus\{0\}$,
the characteristic curve of \eqref{eq:scalar} starting at $\xi$ over
$\gamma(0)$ is the solution of
\begin{equation}\label{eq:flow}
 \dot y(t)=-Q(\gamma(t),y(t))\,\dot\gamma(t),\qquad y(0)=\xi .
\end{equation}
The ODE is meaningful for any such $Q$; the name refers to the case in
which $Q$ solves \eqref{eq:scalar}, where the curve has the two properties
recorded in Lemma~\ref{prop:flow}.

\begin{lemma}\label{prop:flow}
Let $Q\in C^\infty(\Omega\times(\R^n\setminus\{0\}))$ be positively
one-homogeneous. Let
$\gamma:[0,1]\to\Omega$ be piecewise $C^1$, let $\xi\in\R^n\setminus\{0\}$,
and put $C=\gamma([0,1])$ and $L=\int_0^1|\dot\gamma|\dd t$. Then the following hold.
\begin{enumerate}
\item[(i)] Problem~\eqref{eq:flow} has a unique solution
$y:[0,1]\to\R^n\setminus\{0\}$, absolutely continuous and of class $C^1$ on
each interval on which $\gamma$ is $C^1$, and
\begin{equation}\label{eq:flow-bound}
 e^{-ML}|\xi|\le|y(t)|\le e^{ML}|\xi|,\qquad M=\max_{C\times\Sn}|Q| .
\end{equation}
\end{enumerate}
If in addition $Q$ satisfies \eqref{eq:scalar}, then:
\begin{enumerate}
\item[(ii)] $Q(\gamma(t),y(t))=Q(\gamma(0),\xi)$ for $0\le t\le1$, and
\begin{equation}\label{eq:flow-explicit}
 y(t)=\xi-\bigl(\gamma(t)-\gamma(0)\bigr)Q(\gamma(0),\xi).
\end{equation}
In particular the endpoint
\begin{equation}\label{eq:Phi}
 \Theta_\gamma(\xi):=y(1)
 =\xi-\bigl(\gamma(1)-\gamma(0)\bigr)Q(\gamma(0),\xi)
\end{equation}
depends only on $\gamma(0)$ and $\gamma(1)$.
\item[(iii)] $\Theta_\gamma$ is a $C^\infty$ diffeomorphism of
$\R^n\setminus\{0\}$, with $\Theta_\gamma^{-1}=\Theta_{\gamma^-}$ and
$\Theta_{\gamma_2*\gamma_1}=\Theta_{\gamma_2}\circ\Theta_{\gamma_1}$, where
$\gamma^-(t)=\gamma(1-t)$. Moreover
\begin{equation}\label{eq:flow-differential}
 \dd\Theta_\gamma(\xi)
 =\Id-\bigl(\gamma(1)-\gamma(0)\bigr)\otimes\dd_yQ(\gamma(0),\xi),
\end{equation}
and
\begin{equation}\label{eq:flow-det}
 \det\dd\Theta_\gamma(\xi)
 =1-\dd_yQ(\gamma(0),\xi)\bigl(\gamma(1)-\gamma(0)\bigr).
\end{equation}
\end{enumerate}
\end{lemma}

\begin{proof}
(i) Put $V(t,y)=-Q(\gamma(t),y)\dot\gamma(t)$.
On each $C^1$ segment, $V$ is continuous in $t$, smooth in $y$, hence locally Lipschitz in
$y$ uniformly on compact time intervals.
 By the standard
existence and uniqueness theorem for ODEs, \eqref{eq:flow} has a unique
local solution with $y(0)=\xi$, defined on a maximal interval $[0,T)$. On
this interval $y(t)\ne0$, because the solution lies in the slit bundle.

Let $C=\gamma([0,1])$. Since $C\times\Sn$ is compact and $Q$ is continuous
there,
\[
 M=\max_{C\times\Sn}|Q|<\infty,
 \qquad
 |Q(x,y)|\le M|y|\quad\text{for }x\in C,\ y\ne0 .
\]
For $0\le t<T$ we have $\dot y=-Q(\gamma,y)\dot\gamma$, hence
\[
 \left|\frac{\dd}{\dd t}\log|y(t)|\right|
 =\frac{\bigl|Q(\gamma(t),y(t))\bigr|\,
        \bigl|\langle y(t),\dot\gamma(t)\rangle\bigr|}
       {|y(t)|^2}
 \le M|\dot\gamma(t)| .
\]
Integrating over $[0,t]$ gives \eqref{eq:flow-bound}. Thus the curve is
contained in the compact set
\[
 \mathcal K
 =\bigl\{(x,y):x\in C,\ e^{-ML}|\xi|\le|y|\le e^{ML}|\xi|\bigr\}
 \subset\Omega\times(\R^n\setminus\{0\}) .
\]
By the extension theorem for ODEs~\cite[Ch.~II, \S 3]{HartmanODE}, a
solution can be extended as long as its trajectory remains in a compact
subset of the domain. Since the trajectory lies in $\mathcal K$ for
$0\le t\le1$, the maximal interval contains $[0,1]$.

(ii) Along the curve,
\[
 \frac{\dd}{\dd t}Q(\gamma(t),y(t))
 =\dot\gamma^kQ_{x^k}-Q\dot\gamma^kQ_{y^k}
 =\dot\gamma^k\bigl(Q_{x^k}-QQ_{y^k}\bigr)=0 ,
\]
so $Q(\gamma(t),y(t))=Q(\gamma(0),\xi)$. Then \eqref{eq:flow} reads
$\dot y=-Q(\gamma(0),\xi)\dot\gamma$, and integration gives
\eqref{eq:flow-explicit}. The right side of \eqref{eq:flow-explicit}
involves $\gamma$ only through its endpoints, so \eqref{eq:Phi} depends
only on $\gamma(0)$ and $\gamma(1)$.

(iii) By (i), every solution of \eqref{eq:flow} exists on the whole
interval $[0,1]$. Hence the time-one map $\Theta_\gamma$ is defined on
$\R^n\setminus\{0\}$, and it is $C^\infty$ by smooth dependence of
solutions on initial data. If $z(s)=y(1-s)$, then
\[
 z'(s)=-\dot y(1-s)=Q(\gamma(1-s),z(s))\,\dot\gamma(1-s),
\]
which is \eqref{eq:flow} along $\gamma^-$ with $z(0)=y(1)$. Hence
$\Theta_{\gamma^-}\circ\Theta_\gamma=\mathrm{id}$, and $\Theta_\gamma$ is a
diffeomorphism. For a concatenation $\gamma_2*\gamma_1$, the curve is
obtained by following that of $\gamma_1$ and then that of $\gamma_2$, so
uniqueness of solutions gives
$\Theta_{\gamma_2*\gamma_1}=\Theta_{\gamma_2}\circ\Theta_{\gamma_1}$.

Finally, \eqref{eq:Phi} is an identity in $\xi$. Taking its differential
with respect to $\xi$ gives \eqref{eq:flow-differential}; applying
Lemma~\ref{lem:rank-one}  to $\Id-h\otimes\omega$ with
$h=\gamma(1)-\gamma(0)$ and $\omega=\dd_yQ(\gamma(0),\xi)$ gives
\eqref{eq:flow-det}.
\end{proof}

\begin{corollary}\label{cor:Tx}
Let $\Omega$ be a connected open set containing $0$, and let $Q$ be a
$C^\infty$ solution of \eqref{eq:scalar} on
$\Omega\times(\R^n\setminus\{0\})$, positively one-homogeneous in $y$,
with $Q(0,\cdot)=\alpha$. Then for every $x\in\Omega$:
\begin{enumerate}
\item[(i)] for any two piecewise $C^1$ paths
$\gamma_1,\gamma_2:[0,1]\to\Omega$ with $\gamma_i(0)=0$ and
$\gamma_i(1)=x$, the endpoint maps $\Theta_{\gamma_1}$ and
$\Theta_{\gamma_2}$ coincide. Hence there is a well-defined map
\[
 \Theta_x:\R^n\setminus\{0\}\longrightarrow\R^n\setminus\{0\},
 \qquad \Theta_x(\xi)=\xi-x\alpha(\xi),
\]
independent of the choice of path;

\item[(ii)] $\Theta_x$ is a $C^\infty$ diffeomorphism of
$\R^n\setminus\{0\}$ onto itself.
\end{enumerate}
\end{corollary}

\begin{proof}
(i) By Lemma~\ref{prop:flow}(ii), for every $\xi\ne0$,
\[
 \Theta_\gamma(\xi)
 =\xi-\bigl(\gamma(1)-\gamma(0)\bigr)Q(\gamma(0),\xi)
 =\xi-x\,\alpha(\xi),
\]
the last equality because $\gamma(0)=0$, $\gamma(1)=x$ and
$Q(0,\xi)=\alpha(\xi)$. The right side does not involve $\gamma$, so all
endpoint maps agree at every $\xi$.

(ii) Each $\Theta_\gamma=\Theta_x$ is a diffeomorphism by
Lemma~\ref{prop:flow}(iii). For the inverse, let $y\ne0$ and put
$\xi=y+xQ(x,y)$. The implicit equation gives $\alpha(\xi)=Q(x,y)$, hence
$\Theta_x(\xi)=\xi-x\alpha(\xi)=y$. Thus $\Theta_x$ is onto, with the stated
inverse.
\end{proof}

Thus \(\Theta_x\) transports an initial direction at the origin to the direction at
\(x\) reached by the characteristic curve; (i) says that this transport does not depend on the path.

\begin{remark}\label{rem:converse}
Conversely, let $Q\in C^\infty(\Omega\times(\R^n\setminus\{0\}))$ be
positively one-homogeneous and constant along the characteristic curve of
every $C^1$ path. Fix $(x,y)$ and $v$, and take
$\gamma(t)=x+tv$. The curve through $y$ exists by
Lemma~\ref{prop:flow}(i) and satisfies $\dot y(0)=-Q(x,y)v$; hence
\[
 0=\left.\frac{\dd}{\dd t}\right|_{t=0}Q(\gamma(t),y(t))
 =v^k\bigl(Q_{x^k}-QQ_{y^k}\bigr)(x,y) .
\]
As $v$ is arbitrary, $Q$ solves \eqref{eq:scalar}.
\end{remark}

\begin{lemma}\label{lem:unique}
Let $Q_1,Q_2$ be $C^\infty$ solutions of \eqref{eq:scalar} on the same
connected open set $\Omega\ni0$, positively one-homogeneous in $y$, with
$Q_1(0,\cdot)=Q_2(0,\cdot)=\alpha$. Then $Q_1=Q_2$ on
$\Omega\times(\R^n\setminus\{0\})$.
\end{lemma}

\begin{proof}
Fix $x\in\Omega$. By Corollary~\ref{cor:Tx}, the endpoint map of a path
from $0$ to $x$ is $\Theta_x$ for both solutions, and it is a diffeomorphism.
For $i=1,2$ and $\xi\ne0$, Lemma~\ref{prop:flow}(ii) gives
\[
 Q_i\bigl(x,\Theta_x(\xi)\bigr)=Q_i(0,\xi)=\alpha(\xi).
\]
Since $\Theta_x$ is onto, every $y\ne0$ equals $\Theta_x(\xi)$ for some $\xi$, and
then $Q_1(x,y)=\alpha(\xi)=Q_2(x,y)$.
\end{proof}

\subsection{The maximal smooth scalar domain}\label{subsec:maximal}

We now prove Theorem~\ref{thmC}, keeping the notation of \eqref{eq:D}.

\begin{proof}[Proof of Theorem~\ref{thmC}]
First we prove that $\calH_\alpha$ is open and convex. For fixed $u\ne0$
the map $x\mapsto\Delta_\alpha(x,u)=1-\dd\alpha_u(x)$ is affine, so
\[
 f(x)=\min_{|u|=1}\Delta_\alpha(x,u)
\]
is concave. If $L=\sup_{|u|=1}|\dd\alpha_u|$, then
$|f(x)-f(x')|\le L|x-x'|$, so $f$ is continuous; hence
$\calH_\alpha=f^{-1}\bigl((0,\infty)\bigr)$ is open, and it is convex
because $f$ is concave: for $x,x'\in\calH_\alpha$ and $0\le t\le1$,
\[
 f\bigl(tx+(1-t)x'\bigr)\ge tf(x)+(1-t)f(x')>0 .
\]
Moreover $f(0)=1$, so $0\in\calH_\alpha$.

Next we show that any solution of \eqref{eq:scalar} on a connected
domain containing the origin is defined on a subset of $\calH_\alpha$.
Let $Q$ be a solution on a connected open set $\Omega\ni0$, fix
$x\in\Omega$, and let $\gamma$ be a path in $\Omega$ from $0$ to $x$. By
Lemma~\ref{prop:flow}, the endpoint map of $\gamma$ is
$\Theta_x(\xi)=\xi-x\alpha(\xi)$, a diffeomorphism of
$\R^n\setminus\{0\}$; hence $\det\dd\Theta_x(\xi)\ne0$ for $\xi\ne0$. By
Lemma~\ref{lem:rank-one} with $\omega=\dd\alpha_\xi$,
\[
 \det\dd\Theta_x(\xi)=\det\bigl(\Id-x\otimes\dd\alpha_\xi\bigr)
 =\Delta_\alpha(x,\xi).
\]
For fixed $\xi\ne0$ the function $x\mapsto\Delta_\alpha(x,\xi)$ is
continuous and nowhere zero on the connected set $\Omega$, and equals $1$
at $x=0$; it is therefore positive on $\Omega$, and
$\Omega\subset\calH_\alpha$. The same lemma gives
$Q(x,\Theta_x(\xi))=\alpha(\xi)$ for every $\xi\ne0$; as $\Theta_x$ is
onto, this determines $Q(x,\cdot)$.

Conversely, we construct a solution on $\calH_\alpha$. Fix
$x\in\calH_\alpha$ and $y\ne0$, and put
\[
 \delta=\min_{|u|=1}\Delta_\alpha(x,u)>0,\qquad
 h(q)=q-\alpha(y+xq).
\]
At every $q$ with $y+xq\ne0$,
\[
 h'(q)=1-\dd\alpha_{y+xq}(x)=\Delta_\alpha(x,y+xq)\ge\delta .
\]
Since $\dd\alpha$ is zero-homogeneous and bounded on $\Sn$, $\alpha$ is
Lipschitz; hence $h$ is absolutely continuous, and the single exceptional
point $y+xq=0$ does not affect the integral. Therefore
$h(q_2)-h(q_1)\ge\delta(q_2-q_1)$ for $q_2>q_1$, so $h$ is strictly
increasing and $h(q)\to\pm\infty$ as $q\to\pm\infty$. Thus $h$ has
exactly one zero; denote it by $Q(x,y)$. By construction,
$Q(x,y)=\alpha(y+xQ(x,y))$.

The point $\xi=y+xQ(x,y)$ is nonzero: otherwise $Q(x,y)=\alpha(0)=0$, and
then $y=-xQ(x,y)=0$, contrary to $y\ne0$. Hence
$E(x,y,q)=q-\alpha(y+xq)$ is $C^\infty$ near $(x,y,Q(x,y))$ and
$E_q=\Delta_\alpha(x,\xi)\ge\delta>0$ there, so the implicit function
theorem makes $Q$ of class $C^\infty$ near $(x,y)$; as $(x,y)$ was
arbitrary, $Q\in C^\infty(\calH_\alpha\times(\R^n\setminus\{0\}))$. For
$t>0$ the function $q\mapsto q-\alpha(ty+xq)$ vanishes at $tQ(x,y)$, so
uniqueness of the zero gives $Q(x,ty)=tQ(x,y)$; and $Q(0,y)=\alpha(y)$.

It remains to verify that $Q$ solves \eqref{eq:scalar} and to prove
maximality. Put $\xi=y+xQ(x,y)$ and differentiate the implicit equation,
first in $y^k$ and then in $x^k$:
\[
 Q_{y^k}=\frac{\alpha_{\xi^k}(\xi)}{\Delta_\alpha(x,\xi)},\qquad
 Q_{x^k}=\frac{Q(x,y)\,\alpha_{\xi^k}(\xi)}{\Delta_\alpha(x,\xi)}
 =Q\,Q_{y^k}.
\]
Hence $Q$ is a solution on $\calH_\alpha$. The determination
$Q(x,\Theta_x(\xi))=\alpha(\xi)$ of the second paragraph applies to this
solution as well, so any solution on a connected $\Omega\ni0$ agrees with
$Q$ on $\Omega$, and $\Omega\subset\calH_\alpha$ was shown there. This
proves the maximality and the uniqueness.
\end{proof}

\begin{remark}\label{rem:convex-alpha}
No positivity or
convexity of $\alpha$ is assumed in Theorem~\ref{thmC}.
If $\alpha$ is convex, then $\dd\alpha_u(x)\le\alpha(x)$ for every
$u\ne0$, by the supporting-hyperplane inequality at $u$ and Euler's
identity. Hence $f(x)=1-\alpha(x)$ and
$\calH_\alpha=\{\alpha<1\}$, the usual domain obtained in \cite{LZ}.
\end{remark}

\begin{remark} \label{rem:scalar-components}
Maximality in Theorem~\ref{thmC} refers to connected domains
containing $0$. For $\alpha(\xi)=2\xi^1$ and $n\ge2$,
\[
 Q(x,y)=\frac{2y^1}{1-2x^1}
\]
solves \eqref{eq:scalar} on $\{x^1\ne1/2\}$. The component $\{x^1<1/2\}$
is $\calH_\alpha$; on $\{x^1>1/2\}$ one has $\Delta_\alpha<0$, and the
reconstruction below gives $F<0$ for positive $\psi$.
\end{remark}

Theorem~\ref{thmC} describes the maximal domain but not its shape. The
next two results show that the shape is constrained: a complete affine
line in the domain forces the initial data to be linear, and for negative
constant curvature no such line can occur.

\begin{proposition}\label{prop:line}
If $\calH_\alpha$ contains an entire affine line $x_c+\R v$ with $v\ne0$,
then $\alpha$ is a linear functional and $\alpha(v)=0$.
\end{proposition}

\begin{proof}
The inequality $1-\dd\alpha_\xi(x_c+tv)>0$, required for all real $t$,
forces $\dd\alpha_\xi(v)=0$ for every $\xi\ne0$: otherwise the left side is
an affine function of $t$ and changes sign. Hence $\alpha$ is constant
along every segment parallel to $v$ that avoids the origin.

For $z\in v^\perp\setminus\{0\}$ we have $\alpha(z+v)=\alpha(z)$. The
function $f(z)=\alpha(z+v)$ is smooth at $z=0$, and
$f(0)=\alpha(v)=\dd\alpha_v(v)=0$. By positive homogeneity and
differentiability at the origin,
\[
 f(z)=\lim_{t\downarrow0}\frac{f(tz)-f(0)}{t}=\dd f_0(z).
\]
Thus $\alpha$ is linear on $v^\perp$; being constant along $v$ and
continuous, it is linear on $\R^n$.
\end{proof}
For the metric, this gives the following consequence.

\begin{corollary}\label{cor:line}
Suppose a connected coordinate domain of a projectively flat Finsler
metric contains an entire affine line. If $\mathbf{K} =0$, the metric is the
pull-back of a Minkowski norm by the projective transformation
\[
 Z(x)=\frac{x}{1-\ell(x)} .
\]
If $\mathbf{K}=const. <0$, such a domain does not exist.
\end{corollary}
(The proof uses the reconstruction formula of
Subsection~\ref{subsec:reconstruction}.)

\begin{proof}
Translate the initial base point to the origin. For $\mathbf{K} =0$,
Theorem~\ref{thmC} and Proposition~\ref{prop:line} give
$\phi(y)=\ell(y)$ for a linear functional $\ell$. On
$\Delta(x)=1-\ell(x)>0$ define $Z$ as above; its inverse is
$x=Z/(1+\ell(Z))$, and a direct differentiation gives
\[
 \dd Z_x(y)=\frac{(\Id-x\otimes\ell)^{-1}y}{\Delta(x)} .
\]
The reconstruction formula of Subsection~\ref{subsec:reconstruction},
\eqref{eq:flat-F}, together with homogeneity, gives
$F(x,y)=\psi(\dd Z_x(y))$, which is the claimed pull-back.

For $\mathbf{K} =-1$, both $\alpha_+$ and $\alpha_-$ would be linear; then
$2\psi=\alpha_+-\alpha_-$ would be linear, contradicting the Minkowski
property of $\psi$.
\end{proof}

In particular, a projectively flat Finsler metric with $\mathbf{K} =0$ on all of
$\R^n$ is a translation-invariant Minkowski metric in the original
coordinates, since $1-\ell(x)>0$ for every $x$ forces $\ell=0$.

\subsection{Metric reconstruction}\label{subsec:reconstruction}

Theorem~\ref{thmC} solves the scalar equation. We now pass to
the pair $(F,P)$ and to the curvature equations \eqref{eq:system}.

\subsubsection{Zero curvature}

Let $\psi$ be a Minkowski norm and $\phi$ a smooth positively one-homogeneous
function. Applying Theorem~\ref{thmC} with $\alpha=\phi$
gives
\begin{equation}\label{eq:flat-F}
 P(x,y)=\phi(\xi),\qquad \xi=y+xP(x,y),\qquad
 F(x,y)=\frac{\psi(\xi)}{\Delta_\phi(x,\xi)},\qquad x\in\calH_\phi .
\end{equation}
This is the reconstruction formula of \cite{Shen} and
\cite[Theorem~1.1]{Li}, where the verification that $(F,P)$ satisfies
\eqref{eq:system} with $\mathbf K=0$ is carried out; we do not repeat it.

\begin{proposition}\label{prop:flat}
On $\calH_\phi\times(\R^n\setminus\{0\})$ the functions in \eqref{eq:flat-F}
are $C^\infty$ and positively one-homogeneous, $F>0$, and $(F,P)$ solves
\eqref{eq:system} with $\mathbf K=0$. Moreover, if a regular solution with
the same initial data is given on a connected open set $\Omega\ni0$, then
$\Omega\subset\calH_\phi$ and the given solution is the restriction of
\eqref{eq:flat-F}.
\end{proposition}

\begin{proof}
Smoothness and positively one-homogeneity follow from
Theorem~\ref{thmC}, and $F>0$ from $\Delta_\phi>0$ and
$\psi>0$. The equations are those verified in \cite{Shen} and
\cite[Theorem~1.1]{Li}.

For the last assertion, let $(\tilde F,\tilde P)$ be a regular solution
on a connected $\Omega\ni0$ with the same initial data. By
Theorem~\ref{thmC}, $\Omega\subset\calH_\phi$ and $\tilde P$ equals the
projective factor $P$ constructed above. Fix $\xi\ne0$ and
$x\in\Omega$, and let $\gamma$ be a $C^1$ path in $\Omega$ from $0$ to
$x$. Lift it by
\[
 \dot y(t)=-\tilde P(\gamma(t),y(t))\,\dot\gamma(t),\qquad y(0)=\xi .
\]
Since $\mathbf K=0$, the equation \eqref{eq:system} gives
$\tilde P_{x^k}=\tilde P\tilde P_{y^k}$, so $\tilde P$ solves the scalar
equation.
By Lemma~\ref{prop:flow}(ii), and since $\gamma(0)=0$, the point
$\xi=y(t)+\gamma(t)\tilde P(\gamma(t),y(t))$ is independent of $t$; in
particular $\tilde P(\gamma(t),y(t))=\phi(\xi)$, and
$\Delta(t):=1-\dd\phi_\xi(\gamma(t))$ satisfies
$\frac{\dd}{\dd t}\Delta(t)=-\dd\phi_\xi(\dot\gamma(t))$.

Since $\tilde P(x,y)=\phi(\xi)$ with $\xi=y+x\tilde P(x,y)$,
differentiating in $y$ gives
\[
 \bigl(1-\dd\phi_\xi(x)\bigr)\dd_y\tilde P=\dd\phi_\xi,
 \qquad\text{i.e.}\qquad
 \Delta(t)\,\dd_y\tilde P=\dd\phi_\xi .
\]
On the other hand, $F_{x^k}=(PF)_{y^k}$ in \eqref{eq:system} gives
\[
 \frac{\dd}{\dd t}\tilde F(\gamma(t),y(t))
 =\tilde F(\gamma(t),y(t))\,\dd_y\tilde P(\dot\gamma(t)).
\]
Combining the last three identities,
\[
 \frac{\dd}{\dd t}\bigl(\tilde F(\gamma(t),y(t))\,\Delta(t)\bigr)=0 .
\]
Since $\tilde F(0,\xi)\Delta(0)=\psi(\xi)$, we obtain
$\tilde F(x,y)=\psi(\xi)/\Delta_\phi(x,\xi)=F(x,y)$. As $\xi$ and $x$
are arbitrary, $\tilde F=F$ on $\Omega$.
\end{proof}

\subsubsection{Negative curvature}

For $\mathbf K=-1$ the system \eqref{eq:system} is equivalent to two
copies of \eqref{eq:scalar}. With
\[
 Q_\pm=P\pm F,\qquad \alpha_\pm=\phi\pm\psi,
\]
adding and subtracting the two equations of \eqref{eq:system} gives
\begin{equation}\label{eq:split}
 (Q_\pm)_{x^k}=Q_\pm(Q_\pm)_{y^k},
\end{equation}
and conversely $P=\tfrac12(Q_++Q_-)$ and $F=\tfrac12(Q_+-Q_-)$ recover
\eqref{eq:system} from \eqref{eq:split}. This splitting is due to
\cite{Shen} and \cite[Theorem~1.2]{Li}.

Let $\calH=\calH_{\alpha_+}\cap\calH_{\alpha_-}$ and let $Q_\pm$ be the
solutions of Theorem~\ref{thmC} with data $\alpha_\pm$.

\begin{proposition}\label{prop:neg}
On $\calH\times(\R^n\setminus\{0\})$ the functions
$P=\tfrac12(Q_++Q_-)$ and $F=\tfrac12(Q_+-Q_-)$ are $C^\infty$ and positively
one-homogeneous, $F>0$, and $(F,P)$ solves \eqref{eq:system} with
$\mathbf K=-1$. Every regular solution with the same initial data on a
connected open set $\Omega\ni0$ is their restriction.
\end{proposition}

\begin{proof}
Smoothness and one-homogeneity follow from
Theorem~\ref{thmC}; the equations from \eqref{eq:split} by
the computation displayed above.

For positivity, put $h_+(q)=q-\alpha_+(y+xq)$ and
$\xi_-=y+xQ_-(x,y)$. Since $h_+$ is strictly increasing and
$h_+(Q_+)=0$,
\[
 h_+(Q_-)=Q_--\alpha_+(\xi_-)=\alpha_-(\xi_-)-\alpha_+(\xi_-)
 =-2\psi(\xi_-)<0=h_+(Q_+),
\]
so $Q_+>Q_-$, hence $F>0$. Positivity does not give strong convexity;
that is a separate condition, treated in Subsection~\ref{subsec:convex}.

For the last assertion, the functions $\tilde P\pm\tilde F$ of a regular
solution solve the scalar equations on $\Omega$ with data $\alpha_\pm$;
by Theorem~\ref{thmC} they coincide with $Q_\pm$ on $\Omega$, hence
$\tilde F=\tfrac12(Q_+-Q_-)=F$ and $\tilde P=\tfrac12(Q_++Q_-)=P$.
\end{proof}

\subsubsection{The maximal connected regular domain}

For either construction let $\calH$ be the scalar domain and put
\begin{equation}\label{eq:O}
 \calO : =\{x\in\calH:\tfrac12\Hess_y(F^2)(x,y)>0
 \text{ for every }y\ne0\}.
\end{equation}
This set is open: $g=\tfrac12\Hess_y(F^2)$ is zero-homogeneous in $y$, so on
the Euclidean sphere its smallest eigenvalue is continuous, and it is
positive at each point of $\calO$. Also $0\in\calO$, because
$F(0,\cdot)=\psi$ is a Minkowski norm. Let $\calO_0$ be the connected
component of $\calO$ containing $0$.

\begin{corollary}\label{cor:max}
The restriction of the construction to $\calO_0$ is the unique maximal
connected regular projectively flat extension of the given initial data
with curvature $\mathbf K\in\{-1,0\}$ in these affine coordinates. Every
other such connected domain containing $0$ is contained in $\calO_0$.
\end{corollary}

\begin{proof}
Existence follows from \eqref{eq:system} and the definition of $\calO_0$.
Let a regular solution be given on a connected $\Omega\ni0$. By
Theorem~\ref{thmC}, applied to $P$ when $\mathbf K=0$ and to
$P\pm F$ when $\mathbf K=-1$, we have $\Omega\subset\calH$, and the
reconstruction formulas identify its metric function with $F$. Its regular
domain is therefore contained in $\calO$, and connectedness places it in
$\calO_0$.
\end{proof}

\subsection{Strong convexity in the initial data}\label{subsec:convex}

The analysis of $\calO$ through the initial data originates in \cite{LZ}.
There it is proved that in dimension two every connected component of
$\calO$ is convex \cite[Theorem~3.17 and~4.19]{LZ}, and that at a smooth boundary
point of a component the degenerate directions of the fibre Hessian are
tangent to the boundary \cite[Proposition~3.16 and~4.18]{LZ}; the second statement
is what makes $\partial\calO$ computable. Both rest on two devices from
\cite[Lemma~3.13 and~4.17]{LZ}: the radial computation, valid when the base point
lies on the ray of the fibre direction, and the re-basing formula, which
allows the initial data to be taken at an arbitrary point of the domain.

The present section develops this circle of ideas in two directions. It
gives, for every pair $(x,\xi)$, an identity that reduces the strong
convexity of $F$ at $y=\Theta_x(\xi)$ to a condition on $\phi$ and $\psi$
alone, thus removing the restriction $x\parallel y$ of the radial
computation; and it applies this identity to determine $\calO$ explicitly
for Randers data in every dimension (Subsection~\ref{sec:randers}), where the
radial case alone does not suffice. In the plane the two approaches meet:
the tangency statement of Section~\ref{sec:boundary} gives an independent proof
of \cite[Theorem~3.17 and 4.19]{LZ} which assumes neither constant flag curvature
nor the convexity of the component of $\calO$.

\subsubsection{Zero curvature}

Fix $x\in\calH_\phi$ and $\xi\ne0$, and put
\[
 \Delta=\Delta_\phi(x,\xi),\qquad
 A=\Id-x\otimes\dd\phi_\xi,\qquad
 H(\xi)=\frac{\psi(\xi)}{\Delta_\phi(x,\xi)},\qquad
 y=\Theta_x(\xi).
\]
Regard $H$ as a function of $\xi$ with $x$ held fixed. By
Lemma~\ref{lem:rank-one}, $A$ is invertible and $A^{-1}x=x/\Delta$.
Define the symmetric bilinear form
\begin{equation}\label{eq:Q}
 \mathcal Q_{x,\xi}
 =(\Hess_\xi H)(\xi)
  +\frac{(\dd_\xi H)(\xi)(x)}{\Delta_\phi(x,\xi)}\,(\Hess\phi_\xi).
\end{equation}

\begin{proposition} \label{prop:hessian}
With the notation above,
\begin{align}
 A^T(\Hess_yF)A&=\mathcal Q_{x,\xi},\label{eq:pullH}\\
 A^Tg_{x,y}A&=\dd H_\xi\otimes\dd H_\xi+H(\xi)\,\mathcal Q_{x,\xi}.\label{eq:pullg}
\end{align}
Consequently $x\in\calO$ if and only if $\mathcal Q_{x,\xi}$ is positive
definite on $\R^n/\R\xi$ for every $\xi\ne0$.
\end{proposition}

\begin{proof}
All derivatives are taken with respect to $\xi$ at the point $\xi$, with
$x$ fixed.

Since $H=F\circ \Theta_x$ and $\Theta_x(\xi)=\xi-x\phi(\xi)$,
\[
 \dd \Theta_x=A,\qquad \dd^2 \Theta_x(u,v)=-x\,\Hess\phi_\xi(u,v).
\]
Differentiating $H=F\circ \Theta_x$ once gives $\dd H_\xi=\dd_yF\circ A$, that
is,
\[
 A^T\dd_yF=\dd H_\xi .
\]
Differentiating once more,
\begin{align*}
 (\Hess_\xi H)(\xi)(u,v)
 &=(\Hess_yF)(Au,Av)+\dd_yF\cdot\dd^2 \Theta_x(u,v) \\
 &=(\Hess_yF)(Au,Av)-\dd_yF(x)\,\Hess\phi_\xi(u,v),
\end{align*}
hence
\begin{equation}\label{eq:pointwise}
 (\Hess_yF)(Au,Av)
 =(\Hess_\xi H)(\xi)(u,v)+\dd_yF(x)\,\Hess\phi_\xi(u,v).
\end{equation}
Because $A^{-1}x=x/\Delta$,
\[
 \dd_yF(x)=\dd_yF\Bigl(A\frac{x}{\Delta}\Bigr)
 =\frac{(A^T\dd_yF)(x)}{\Delta}
 =\frac{\dd H_\xi(x)}{\Delta}.
\]
Substituting this into \eqref{eq:pointwise} and comparing with the
definition \eqref{eq:Q} of $\mathcal Q_{x,\xi}$,
\[
 (\Hess_yF)(Au,Av)=\mathcal Q_{x,\xi}(u,v)
 \qquad\text{for all }u,v\in\R^n .
\]
Since $(A^T(\Hess_yF)A)(u,v)=(\Hess_yF)(Au,Av)$ by the definition of the
pull-back, this is the identity \eqref{eq:pullH} of bilinear forms.

For \eqref{eq:pullg}, Lemma~\ref{lem:angular} gives
$g_{x,y}=\dd F\otimes\dd F+F\Hess F$. Pulling back by $A$ and using
$F(y)=H(\xi)$,
\[
 A^Tg_{x,y}A
 =(A^T\dd_yF)\otimes(A^T\dd_yF)+H(\xi)\,A^T(\Hess_yF)A
 =\dd H_\xi\otimes\dd H_\xi+H(\xi)\,\mathcal Q_{x,\xi}.
\]

For the criterion, note that a $C^2$ positively one-homogeneous function $f$
satisfies $\Hess f(y,\cdot)\equiv0$, since Euler's identity
$f_{y^i}(y)y^i=f(y)$ gives $f_{y^iy^j}(y)y^i=0$. Hence $\Hess_yF$ and
$\mathcal Q_{x,\xi}$ annihilate $y$ and $\xi$, respectively, and define
forms on $\R^n/\R y$ and $\R^n/\R\xi$. As $A\xi=y$ and $A$ is invertible,
$A$ identifies these quotients, and \eqref{eq:pullH} gives
$\mathcal Q_{x,\xi}(u,u)=(\Hess_yF)(Au,Au)$. So positivity of
$\mathcal Q_{x,\xi}$ on $\R^n/\R\xi$ is equivalent to positivity of
$\Hess_yF$ on $\R^n/\R y$, hence, by Lemma~\ref{lem:angular}, to
$g_{x,y}>0$, i.e., to $x\in\calO$.
\end{proof}

Expanding \eqref{eq:Q} is useful in computations. Write
$C=\Hess\phi_\xi$ and $T=\dd^3\phi_\xi$. Since
$\Delta=1-\sum_l\phi_{\xi^l}(\xi)x^l$,
\[
 \dd\Delta_\xi=-C(x,\cdot),\qquad \Hess\Delta_\xi=-T(x,\cdot,\cdot),
\]
and differentiating $H=\psi/\Delta$ once and then once more gives
\[
 \dd H_\xi=\frac{\dd\psi_\xi}{\Delta}
 +\frac{\psi(\xi)C(x,\cdot)}{\Delta^2},
\]
\begin{align*}
 (\Hess_\xi H)(\xi)
 &=\frac{\Hess\psi_\xi}{\Delta}
 +\frac{\dd\psi_\xi\otimes C(x,\cdot)
        +C(x,\cdot)\otimes\dd\psi_\xi
        +\psi(\xi)T(x,\cdot,\cdot)}{\Delta^2} \\
 &\quad +\frac{2\psi(\xi)C(x,\cdot)\otimes C(x,\cdot)}{\Delta^3}.
\end{align*}
Substituting these and $\dd H_\xi(x)$ into \eqref{eq:Q} yields
\begin{align}\label{eq:expandedQ}
 \mathcal Q_{x,\xi}(u,v)={}&\frac{\Hess\psi(u,v)}\Delta\notag\\
 &+\frac{\dd\psi(u)b(v)+\dd\psi(v)b(u)+\psi T(x,u,v)+\dd\psi(x)C(u,v)}{\Delta^2}\notag\\
 &+\frac{2\psi b(u)b(v)+\psi C(x,x)C(u,v)}{\Delta^3},
\end{align}
where $b(u)=C(x,u)$. Thus third derivatives of $\phi$ enter the criterion,
while $\Delta\ne0$ alone controls only the smoothness of $F$.

At a boundary point $x_*\in\partial\calO_0\cap\calH$ the smallest angular
eigenvalue vanishes in at least one direction: it cannot be negative by
continuity from $\calO_0$, and cannot be positive everywhere by openness of
$\calO$. Strong convexity therefore degenerates at such a point, and the
metric does not continue regularly through it.

\subsubsection{Negative curvature}\label{subsubsec:negative}

For $\mathbf K=-1$ the metric function is a difference of two scalar
solutions, and its fibre Hessian is computed from theirs by
Lemma~\ref{lem:impHess}.

Let now
\[
 \alpha_\pm=\phi\pm\psi,
\]
and let $Q_\pm$ be the solutions of Theorem~\ref{thmC} with
data $\alpha_\pm$. By Proposition~\ref{prop:neg},
\[
 F=\tfrac12(Q_+-Q_-)
\]
is positive on $\calH=\calH_{\alpha_+}\cap\calH_{\alpha_-}$ and solves
\eqref{eq:system} with $\mathbf K=-1$.

Note that the two Hessians enter with opposite signs. Positivity of
$\Hess_yQ_+$ and of $\Hess_yQ_-$ separately is neither assumed nor
sufficient, and there is no analogue of the single pull-back
\eqref{eq:pullH}: the maps $A_+$ and $A_-$ do carry $\R^n/\R\xi_+$ and
$\R^n/\R\xi_-$ to $\R^n/\R y$, but the two forms are evaluated at
different points $\xi_+\ne\xi_-$ and are pulled back by the different
maps $A_+^{-1}$ and $A_-^{-1}$. They must be compared on $\R^n/\R y$
after these two pull-backs, not by identifying $\xi_+$ with $\xi_-$.

\begin{proposition}\label{prop:negative}
With the notation above, put
\[
 \xi_\pm=y+xQ_\pm(x,y),\qquad
 \Delta_\pm=\Delta_{\alpha_\pm}(x,\xi_\pm),\qquad
 A_\pm=\Id-x\otimes\dd\alpha_\pm|_{\xi_\pm}.
\]
Then
\begin{equation}\label{eq:negativeH}
 \Hess_yF=\frac12\left(
 \frac{A_+^{-T}(\Hess\alpha_+)_{\xi_+}A_+^{-1}}{\Delta_+}
 -\frac{A_-^{-T}(\Hess\alpha_-)_{\xi_-}A_-^{-1}}{\Delta_-}
 \right).
\end{equation}
Since $F>0$ on $\calH$, Lemma~\ref{lem:angular} and
\eqref{eq:negativeH} show that
\[
 x\in\calO
 \iff \eqref{eq:negativeH}\ \text{is positive definite on}\ \R^n/\R y
 \ \text{for every }y\ne0 .
\]
\end{proposition}

\begin{proof}
Apply Lemma~\ref{lem:impHess} to $\alpha_+$ and to $\alpha_-$. Since
$\Hess_yF=\frac12(\Hess_yQ_+-\Hess_yQ_-)$, this gives
\eqref{eq:negativeH}.

Since $F>0$ on $\calH$, Lemma~\ref{lem:angular} and the identity
\eqref{eq:negativeH} show that $x\in\calO$ precisely when the form
\eqref{eq:negativeH} is positive definite on $\R^n/\R y$ for every
$y\ne0$.
\end{proof}

For $\mathbf{K} \le0$, Theorem~\ref{thmC} thus gives the exact
maximal connected scalar domain containing the initial point.
The angular tests above determine the strong-convexity locus within
it. Positive curvature requires a separate continuation argument,
since the holomorphic generator and its implicit branch must also
be controlled.

\section{Continuation in positive curvature}\label{sec:continuation}

For $\mathbf K\le0$ Theorem~\ref{thmC} determined the maximal
connected base domain of the scalar equation, and the set $\calO$ then
located the subset on which the reconstructed function is a Finsler
metric. No analogue is available for $\mathbf K=1$. There the scalar
equation is replaced by
\[
 Q=h(y+xQ),\qquad Q=P+iF,
\]
where $h$ is the holomorphic continuation of $\phi+i\psi$. Theorem~\ref{thmA}
guarantees that $h$ exists near the real fibre, but its domain of
holomorphy $W\subset\C^n$ need not be all of $\C^n$, and three things can
fail at a finite distance: the branch $z=y+xQ$ may reach $\partial W$, the
denominator $\Delta=1-\dd h_z(x)$ may vanish, or the reconstructed metric
may lose strong convexity.
Whenever one of them fails, the branch cannot be continued in this chart; a larger chart may still be possible.

In a fixed holomorphic chart, the criterion below shows that these are the only obstructions along a finite arc.
 If, on the tail of the arc, the branch stays in a compact
subset of $W$, the denominator stays bounded away from zero, $Q$ stays
bounded with $\Im Q$ bounded below, and the fundamental tensor stays
uniformly positive, then the solution extends to a neighbourhood of the
endpoint, and the extension is real analytic. The criterion is conditional
by design: it does not assert that the bounds hold, and it therefore does
not produce a maximal domain for $\mathbf K=1$ comparable to
Theorem~\ref{thmC}.

Since $h$ is holomorphic, Lemma~\ref{lem:impHess} applies with
$\alpha=h$; the Hessian is complex bilinear, and
\begin{equation}\label{eq:complexH}
 \Hess_yQ=\frac{A^{-T}(\Hess h)_zA^{-1}}{\Delta},
 \qquad z=y+xQ,\quad \Delta=1-\dd h_z(x).
\end{equation}
The formula is used in the proof of Theorem~\ref{thm:continue}: its
right-hand side is continuous in $(x,u,Q)$, so along the path
$\Hess_yQ(x(t),u)\to\Hess_yQ(x_*,u)$, and the fundamental tensors
converge. This is how strong convexity passes to the limit.

\begin{theorem}\label{thm:continue}
Let $F$ be a regular projectively flat Finsler metric on a domain
$\Omega\subset\R^n$ with $\mathbf K=1$, and let
$P$ be its projective factor. Let $x:[0,1)\to\Omega$ be a $C^1$ path of
finite Euclidean length with $x(t)\to x_*$ as $t\to1$.

Suppose $Q=P+iF$ satisfies $Q=h(y+xQ)$ near the path tail in every
nonzero fibre direction. Here $h$ is one holomorphic function on an
open set $W\subset\C^n$ and satisfies~\eqref{eq:euler-h}. Define
\[
 z(t,u)=u+x(t)Q(x(t),u),\qquad
 \Delta(t,u)=1-\dd h_{z(t,u)}(x(t)).
\]
Assume there are constants $M,m,\delta,\kappa>0$ and a compact set
$C\Subset W$ such that, for every $t$ sufficiently close to $1$ and
every $|u|=1$,
\begin{align}
 |Q(x(t),u)|\le M, &\qquad\qquad \im Q(x(t),u) \ge m,\label{eq:cont-bounds}\\
 z(t,u)\in C,   &\qquad\qquad |\Delta(t,u)| \ge\delta,\label{eq:cont-compact}\\
 F(x(t),u)\Hess_y F(x(t),u)(v,v)\ge\kappa |v|^2
 & \qquad\text{for every }v\in\R^n\text{ with }\langle v,u\rangle=0.
 \label{eq:cont-convex}
\end{align}
Then this branch has a regular real-analytic extension to a neighbourhood of $x_*$. It agrees with the original solution on a neighbourhood of the path tail.
\end{theorem}

\begin{proof}
Since $C\Subset W$, there is a compact set $K$ with $C\subset K\subset W$.
All derivatives of $h$ are bounded on $K$; in particular, $|\dd h|\le L$
there. From~\eqref{eq:firstderivatives}, along the path,
\[
 \left|\frac{\dd}{\dd t}Q(x(t),u)\right|
 =\left|\frac{Q\,\dd h_{z(t,u)}(\dot x(t))}{\Delta(t,u)}\right|
 \le \frac{ML}{\delta}|\dot x(t)|.
\]
Since the path has finite length, $Q(x(t),u)$ converges, uniformly for $u\in\Sn$, to a limit $Q_*(u)$, which is continuous. Put
\[
 z_*(u)=u+x_*Q_*(u),\qquad
 \Delta_*(u)=1-\dd h_{z_*(u)}(x_*).
\]
Then $z_*(u)\in C$, and the limiting equation and bounds are
\[
 Q_*(u)=h(z_*(u)),\qquad |\Delta_*(u)|\ge\delta,\qquad \im Q_*(u)\ge m.
\]
Set $E(x,y,q)=q-h(y+xq)$. At each $(x_*,u,Q_*(u))$,
\[
 E(x_*,u,Q_*(u))=0,\qquad
 |E_q(x_*,u,Q_*(u))|=|\Delta_*(u)|\ge\delta.
\]
The holomorphic implicit function theorem gives a local root near each
point of this compact graph.

These roots have a uniqueness neighbourhood of uniform size. Indeed,
near a point of the graph choose a neighbourhood of the parameters and
a disc in the $q$ variable on which $E_q$ differs from its value at
that point by less than $\delta/2$. If $q_1,q_2$ are two roots in this
disc for the same parameters, then
\[
 0=(q_1-q_2)\int_0^1
 E_q(x,y,q_2+s(q_1-q_2))\,\dd s.
\]
Their joining segment stays in the disc. The integral differs from
the derivative at the graph point by less
than $\delta/2$. Its modulus is therefore at least $\delta/2$,
which implies $q_1=q_2$.
Compactness makes the neighbourhood size uniform along the graph.

Continuity of $Q_*$ identifies each local root with the prescribed root
at the endpoint $x=x_*$. A finite cover of the sphere, followed by shrinking the
neighbourhoods, makes the roots agree on overlaps and gives a common
base neighbourhood.

Euler's identity for $h$ and local uniqueness give one-homogeneity of the patched branch, exactly as in Theorem~\ref{thm:data}; extend it radially to all $y\ne0$. First differentiation gives~\eqref{eq:burgers}.

By~\eqref{eq:firstderivatives} and~\eqref{eq:complexH}, \(Q\), its
first two fibre derivatives, and hence \(F\) and \(\Hess_yF\), depend
continuously on \((x,u,Q)\). Therefore
\(F(x(t),u)\Hess_yF(x(t),u)\) converges to
\(F(x_*,u)\Hess_yF(x_*,u)\), and the inequality in
\eqref{eq:cont-convex} passes to the limit. By
Lemma~\ref{lem:angular}, \(g_{x_*,u}>0\) for every \(u\ne0\);
continuity and compactness of \(\Sn\) then give a uniform positive
lower bound on a smaller base neighbourhood. Taking \(P=\re Q\) and
\(F=\im Q\) gives a regular solution of~\eqref{eq:system}.

Uniform convergence places the original roots in these same implicit neighbourhoods for all sufficiently large $t$. Implicit uniqueness identifies the branches near each point of that tail, not merely at the endpoint.
\end{proof}

\begin{remark}\label{rem:noautomatic}
The theorem assumes a single holomorphic generator on a neighbourhood of the compact complex image. More general continuation may require several holomorphic charts and may have monodromy. The bounds on $\Delta$ and on the angular eigenvalues are also hypotheses. Extension to a projective compactification would require control in the charts meeting infinity. In particular, the closed-hemisphere regularity hypothesis in \cite[Theorem~5.6]{LZ} remains a separate condition.
\end{remark}

\section{Regularity boundaries and Randers initial data}
\label{sec:boundary}

In Section~\ref{sec:real} we reduced strong convexity to positivity of
the angular Hessian and introduced the regularity locus \(\calO\).
This section studies the boundary of \(\calO\) using only Hamel's
equation; the flag curvature need not be constant. The first
subsection proves that the fibre Hessian is invariant along straight
lines, and derives planar convexity and boundary tangency from this
invariance. The second computes \(\calO\) for equal Randers data in
every dimension. 

\subsection{Geometry of the regularity boundary}
\label{subsec:boundary-geometry}

In the constant flag curvature setting, Li and Zhao \cite[Theorems~3.17 and~4.19]{LZ} established the convexity of regularity components for nonpositive curvature in the plane. A natural question is whether this geometric property is a consequence of the constant curvature constraint, or merely of projective flatness itself. We show in this subsection that it is the latter. By relying solely on Hamel's equation, we prove that the fibre Hessian is invariant along straight lines (Theorem~\ref{thm:line-invariance}).
This identity yields planar convexity of the regularity components for
smooth nonnegative Hamel solutions on a convex ambient domain, without a
constant-curvature assumption (Corollary~\ref{cor:planar}). In dimensions \(n \ge 3\), this invariance manifests as the tangency of degenerate directions at smooth regularity boundaries (Theorem~\ref{thm:boundary}).

\begin{theorem}\label{thm:line-invariance}
Let $F=F(x,y)\geq 0$ and $P=P(x,y) $ be smooth, positively one-homogeneous functions on the slit bundle over an open set $\Omega$, and suppose $F_{x^k}=(PF)_{y^k}$. Strong convexity is not assumed. Then the matrix-valued function $t\mapsto\Hess_yF(x+ty,y)$ is constant on
each interval of the line $\{x+ty:t\in\R\}$ that is contained in $\Omega$.
\end{theorem}

\begin{proof}
Recall that $\Dop=y^k\partial_{x^k}$. Euler's formula for the two-homogeneous function $PF$ gives
\[
 \Dop F=y^k(PF)_{y^k}=2PF.
\]
Differentiate twice in $y$, accounting for the $y$ in $\Dop$:
\[
 \partial_{y^i}\partial_{y^j}(\Dop F)
 =\Dop(F_{y^iy^j})+F_{x^iy^j}+F_{x^jy^i}.
\]
Each of the last two terms is $(PF)_{y^iy^j}$ (they are equal because
$PF$ is $C^2$). Comparing with
$\partial_{y^i}\partial_{y^j}(2PF)=2(PF)_{y^iy^j}$ gives
$\mathcal D(F_{y^iy^j})=0$.
Since $\mathcal D$ is the derivative along $t\mapsto x+ty$, each entry
$F_{y^iy^j}$ is constant on every interval of this line contained in
$\Omega$; hence so is $\Hess_yF$.
\end{proof}

For \(p\in\Omega\) and \(y\ne0\), \(p+\R y\) is the affine line through
\(p\) in direction \(y\). By Theorem~\ref{thm:line-invariance}, the
Hessian \(\Hess_yF(\cdot,y)\) is constant along this line.

\begin{theorem}\label{thm:boundary}
Assume that $F$ and $P$ satisfy the hypotheses of
Theorem~\ref{thm:line-invariance} on a convex open domain $\Omega$.
Put
\[
 \calO=\{x\in\Omega:g_{x,y}>0\ \text{for all }y\ne0\}.
\]
Assume in addition that $\calO\ne\varnothing$, and let $U$ be a connected component of $\calO$.  Then
$F(x,y)>0$ for every $x\in\Omega$ and $y\ne0$.
At every $p\in\partial U\cap\Omega$, there is a direction $y\ne0$ for which $\Hess_yF(p,y)$ has a null vector outside $\R y$. For every such direction,
\begin{equation}\label{eq:excluded-line}
 (p+\R y)\cap\Omega\cap\calO=\varnothing.
\end{equation}
If, near $p$, $\partial U$ is a $C^1$ hypersurface with $U$ on one side, then every such $y$ is tangent to $\partial U$ at $p$. No convexity of $U$ is assumed.
\end{theorem}

\begin{proof}
Fix $x_1\in\calO$ and $\xi_1\ne0$, so
$F(x_1,\xi_1)>0$. Given $x\in\Omega$ and $y\ne0$, let $\gamma$ be a
$C^1$ path in $\Omega$ from $x_1$ to $x$, and solve
\[
 \dot y(t)=-P(\gamma(t),y(t))\,\dot\gamma(t),\qquad y(1)=y,
\]
backward on $[0,1]$. The homogeneity estimate used in the proof of
Lemma~\ref{prop:flow}(i)---it uses only smoothness and positive
one-homogeneity of $P$, not the scalar equation---gives a solution with
$y(t)\ne0$. Along this lift,
\[
 \frac{\dd}{\dd t}F(\gamma(t),y(t))
 =F(\gamma(t),y(t))\,\dd_yP(\dot\gamma(t)),
\]
because $F_{x^k}=(PF)_{y^k}$. Hence
\[
 F(x,y)=F(x_1,y(0))\,
 \exp\Bigl(\int_0^1\dd_yP(\dot\gamma(t))\,\dd t\Bigr)>0 .
\]

Since \(p\in\partial U\) and \(U\) is open, choose \(x_k\in U\) with
\(x_k\to p\). For fixed \(y\ne0\), the forms
\(\Hess_yF(x_k,y)\) are positive definite on \(\R^n/\R y\) and
converge to \(\Hess_yF(p,y)\). Hence \(\Hess_yF(p,y)\) is positive
semidefinite on \(\R^n/\R y\). It vanishes on the radial line by
homogeneity.

For \(|y|=1\), let
\[
\lambda(x,y)=\min_{\substack{v\perp y\\ |v|=1}}
\Hess_yF(x,y)(v,v).
\]
This is continuous on \(\Omega\times\Sn\). If \(\lambda(p,y)>0\) for
every \(y\), then \(\min_{|y|=1}\lambda(p,y)>0\) by compactness, and
continuity gives \(\lambda>0\) near \(p\). Then \(g>0\) near \(p\), so
\(p\in\calO\), a contradiction. Hence \(\lambda(p,y)=0\) for some
\(y\ne0\), i.e.\ \(\Hess_yF(p,y)\) has a nonradial null vector.

Fix such a \(y\). For \(p+t y\in\Omega\),
Theorem~\ref{thm:line-invariance} gives
\(\Hess_yF(p+t y,y)=\Hess_yF(p,y)\). The null vector persists, so
\(p+t y\notin\calO\). Hence
\((p+\R y)\cap\Omega\cap\calO=\varnothing\).

If \(\partial U\) is \(C^1\) near \(p\) and \(y\) is not tangent to it,
the line \(p+\R y\) crosses the boundary and contains points of \(U\)
near \(p\), contradicting the exclusion. Thus \(y\) is tangent to $\partial U$ at $p$.
\end{proof}

\begin{corollary}\label{cor:planar}
Under the hypotheses of Theorem~\ref{thm:boundary}, if $n=2$, every connected component of $\calO$ is convex.
\end{corollary}
\begin{proof}
Let \(a,b\in U\), and suppose the segment \([a,b]\) is not contained in
\(U\). Since \(U\) is open, there is \(p\in[a,b]\cap\partial U\); since
\(\Omega\) is convex, \(p\in\Omega\).
At \(p\), choose a degenerate direction \(y\) as in
Theorem~\ref{thm:boundary}, and set \(L=p+\R y\).
Since \(\Omega\) is convex, \(L\cap\Omega\) is an open segment or an open line containing
\(p\). For every \(q\in L\cap\Omega\),
Theorem~\ref{thm:line-invariance} gives
\[
\Hess_yF(q,y)=\Hess_yF(p,y).
\]
Hence the nonradial null vector of \(\Hess_yF(p,y)\) is also a null
vector of \(\Hess_yF(q,y)\). By Lemma~\ref{lem:angular}, \(g_{q,y}\)
is not positive definite, so \(q\notin\calO\). Therefore
\(L\cap\Omega\cap\calO=\varnothing\), and since \(U\subset\calO\),
\(L\cap U=\varnothing\).

 As \(U\) is connected, it lies in one of the two open
half-planes bounded by this line. Both \(a\) and \(b\) lie in that
half-plane, so \([a,b]\) does too, contradicting \(p\in[a,b]\) and
\(p\in p+\R y\).
\end{proof}

Corollary~\ref{cor:planar} proves convexity of the regularity components
only in the plane; its proof uses the separation of the plane by a line,
which fails in higher dimensions. The next example shows that convexity
does indeed fail there. It is a Hilbert metric whose boundary has two
flat points: the angular Hessian of the boundary function degenerates
along a single axis, and the Hessian formula of
Proposition~\ref{prop:negative} identifies the exact regularity locus.

\begin{proposition}
\label{prop:nonconvex}
For every \(n\ge3\) there are real-analytic, strongly convex
one-fibre data for a projectively flat metric of flag curvature
\(-1\) whose maximal connected regular domain in the affine chart is
nonconvex.
\end{proposition}

\begin{proof}
Fix a Euclidean unit vector \(e\) and set
\[
 f(y)=|y|+\frac{\langle e,y\rangle^2}{|y|},
 \qquad
 H=\{x\in\R^n:f(x)<1\},
\]
with \(f(0)=0\). For \(|u|=1\) and \(v\perp u\),
\begin{equation}\label{eq:nonconvex-angular}
 \Hess f_u(v,v)
 =(1-\langle e,u\rangle^2)|v|^2+2\langle e,v\rangle^2 .
\end{equation}
Thus $\Hess f_u\ge0$ for every $u\ne0$, so $f$ is convex on
$\R^n\setminus\{0\}$; being even, positively one-homogeneous and
continuous at $0$, it is convex on all of $\R^n$, and $H$ is open,
bounded and convex.

The form \eqref{eq:nonconvex-angular} degenerates exactly on $\R e$:
for $\xi\in\R e\setminus\{0\}$ it vanishes on $\xi^\perp$.
Since $f\ge0$ and $f$ is convex, $f^2$ is convex; it is not strongly
convex, because $\Hess(f^2)=2f\,\Hess f+2\,\dd f\otimes\dd f$ vanishes
on $\xi^\perp$ when $\xi\in\R e\setminus\{0\}$. Hence $f$ is not a
Minkowski norm.
Since the reconstruction below satisfies $F(0,\cdot)=f$, the origin is
not a regular point. This is why, at the end of the proof, the initial
data are taken at a point $x_0\notin\R e$ rather than at the origin.

Apply Theorem~\ref{thmC} to
\[
\alpha_+=f,\qquad \alpha_-=-f.
\]
Convexity and Euler's identity give
\(\sup_{|u|=1}\dd f_u(x)=f(x)\); since \(f\) is even,
\(\inf_{|u|=1}\dd f_u(x)=-f(x)\). Hence
\[
\calH_{\alpha_+}=\calH_{\alpha_-}=H.
\]
Let \(Q_\pm\) be the corresponding solutions of \eqref{eqQ_thmC} for \(\alpha = \alpha_{\pm}\). At \(x=0\),
\(Q_+=f(y)>0\) and \(Q_-=-f(y)<0\); by continuity and the absence of
zeros, \(Q_+>0>Q_-\) on \(H\). Put
\[
F=\tfrac12(Q_+-Q_-),\qquad P=\tfrac12(Q_++Q_-).
\]
Then \(F>0\), and \((F,P)\) is a real-analytic solution of
\eqref{eq:system} with \(\mathbf K=-1\).  Geometrically, the two points \(x+y/Q_\pm\) lie on \(\partial H\), so \(F\) is the Hilbert metric of \(H\) wherever it is strongly convex.

Write
\[
 \xi_\pm=y+xQ_\pm,\qquad
 \Delta_\pm=\Delta_{\alpha_\pm}(x,\xi_\pm),\qquad
 A_\pm=\Id-x\otimes\dd\alpha_{\pm,\xi_\pm}.
\]
Since \(x\in H=\calH_{\alpha_\pm}\), we have \(\Delta_\pm>0\), so
\(A_\pm\) is invertible. By \eqref{eq:negativeH},
\begin{equation}\label{eq:nonconvex-Hess}
 2\Hess_yF=
 \frac{A_+^{-T}(\Hess f)_{\xi_+}A_+^{-1}}{\Delta_+}
 +\frac{A_-^{-T}(\Hess f)_{\xi_-}A_-^{-1}}{\Delta_-}.
\end{equation}
The two signs add because \(\alpha_-=-f\).

We record the angular behaviour of \(\Hess f\). For \(\xi\ne0\),
\(\Hess f_\xi\) is positive semidefinite and annihilates \(\xi\):
writing \(v=c\xi+w\) with \(w\perp\xi\),
\[
 \Hess f_\xi(v,v)=\Hess f_\xi(w,w) .
\]
By \eqref{eq:nonconvex-angular}, the induced form on \(\R^n/\R\xi\) is
positive definite if \(\xi\notin\R e\). If \(\xi\in\R e\), then
\(u=\pm e\), and the form vanishes for \(w\perp e\); the induced form
is only semidefinite.

Put
\[
 T_\pm=
 \frac{A_\pm^{-T}(\Hess f)_{\xi_\pm}A_\pm^{-1}}{\Delta_\pm}.
\]
Each \(T_\pm\) is positive semidefinite. Euler's identity gives
\[
 A_\pm\xi_\pm
 =\xi_\pm-x\,\dd\alpha_{\pm,\xi_\pm}(\xi_\pm)
 =\xi_\pm-xQ_\pm=y,
\]
so \(A_\pm\) identifies \(\R^n/\R\xi_\pm\) with \(\R^n/\R y\).
Therefore \(T_\pm\) is positive definite on \(\R^n/\R y\) if and only
if \(\xi_\pm\notin\R e\); if \(\xi_\pm\in\R e\), then \(T_\pm\) is not
positive definite on the quotient.

Now let \(x\notin\R e\). If both \(\xi_+,\xi_-\in\R e\), then
\[
 \xi_+-\xi_-=x(Q_+-Q_-)=2xF\in\R e.
\]
Since \(F>0\), this forces \(x\in\R e\), a contradiction. Hence at
least one of \(\xi_\pm\) lies outside \(\R e\). The corresponding
\(T_\pm\) is positive definite on \(\R^n/\R y\), while the other is
positive semidefinite; their sum is positive definite. Thus
\(x\in\calO\).

Conversely, let \(x\in H\cap\R e\), say \(x=\lambda e\), and take
\(y=e\). The equations for \(Q_\pm\) and the data \(\alpha_\pm=\pm f\)
are invariant under every rotation \(R\) with \(Re=e\); the domain
\(H\) is invariant as well. By uniqueness,
\(Q_\pm(Rx,Ry)=Q_\pm(x,y)\). Hence
\[
 R\xi_\pm(x,y)
 =Ry+Rx\,Q_\pm(Rx,Ry)
 =y+xQ_\pm(x,y)=\xi_\pm(x,y).
\]
For \(x=\lambda e\) and \(y=e\), this says that \(\xi_\pm\) is fixed
by every rotation with \(Re=e\); therefore \(\xi_\pm\in\R e\).
Since $\xi_\pm\in\R e$, we have $\dd\alpha_{\pm,\xi_\pm}(v)=0$ and
$(\Hess f)_{\xi_\pm}(v,\cdot)=0$ for every $v\perp e$; hence $A_\pm v=v$
and $T_\pm(v,v)=0$. Thus $(T_++T_-)(v,v)=0$ for every $v\perp e$, and
$v\perp e$ is nonzero in $\R^n/\R e$. Therefore $\Hess_yF$ is not
positive definite on the angular quotient, and $x\notin\calO$.

Consequently
\begin{equation}\label{eq:nonconvex-O}
 \calO=H\setminus\R e.
\end{equation}
In particular, although \(f\) itself is only convex, the fibre at every
\(x_0\in\calO\) is strongly convex.

Since $n\geq 3$ and $\R e$ has codimension $n-1 \geq 2$, it does not disconnect the open convex
set $H$; hence $\calO=H\setminus\R e$ is path connected. Moreover, $\calO$ is not
convex: for $v\perp e$, $v\ne0$, and small $t>0$ the points $\pm tv$ lie
in $\calO$, while the segment joining them contains $0\in\R e$.

Finally, choose \(x_0\in\calO\). The one-fibre data at \(x_0\) are
\[
 \psi=F(x_0,\cdot),\qquad \phi=P(x_0,\cdot).
\]
Since \(x_0\in\calO\), \(\psi\) is a Minkowski norm and \(\phi\) is
smooth and positively one-homogeneous; both are real analytic.
For either branch $\alpha=\alpha_\pm$, let $\alpha'(y)=Q(x_0,y)$ be the
translated data, write $\xi=y+xQ(x,y)$, $Q=Q_\pm$, and put
$\eta=\xi-x_0\alpha(\xi)$, $x=x_0+x'$.
Differentiating the
identity \(Q(x_0,\eta)=\alpha(\xi)\) implicitly gives
\[
 1-\dd_yQ(x_0,\eta)(x')
 =\frac{1-\dd\alpha_\xi(x)}{1-\dd\alpha_\xi(x_0)} .
\]
Since \(x_0\in\calH_{\alpha_\pm}\), the denominator is positive.
As $\eta$ runs over $\R^n\setminus\{0\}$, so does $\xi$, since
$\xi\mapsto\xi-x_0\alpha(\xi)$ is a diffeomorphism of
$\R^n\setminus\{0\}$; hence $\calH_{\alpha'}=H-x_0$, and its regular locus is \(\calO-x_0\). Since \(\calO\) is path connected, \(\calO-x_0\) is connected; it is the
component of the regular locus that contains the translated origin
(the image of \(x_0\)). By Corollary~\ref{cor:max}, it
is the maximal connected regular extension of the translated data. It is
nonconvex because \(\calO\) is.
\end{proof}

\begin{remark}\label{rem:planar-contrast}
For $n=2$ the same construction gives the two convex components
$H\cap\{\langle e_\perp,x\rangle>0\}$ and
$H\cap\{\langle e_\perp,x\rangle<0\}$, where $e_\perp\perp e$, in
accordance with Corollary~\ref{cor:planar}.
\end{remark}

\subsection{Application: Randers initial data}\label{sec:randers}

The calculation of the exact regularity domain for equal Randers initial data in the plane (\(n=2\)) was first carried out by Li and Zhao \cite[Proposition~3.12]{LZ}. In this subsection, we extend their result to arbitrary dimensions \(n \ge 2\). The passage to higher dimensions requires two new ingredients. First, the angular space has dimension \(n-1\), so checking a single nonradial direction (as in the planar case) is no longer sufficient; our proof requires a spectral analysis of a rank-two symmetric operator to test all angular directions. Second, the planar convexity argument, which relies on a line separating the plane, does not carry over to \(n \ge 3\). We establish the higher-dimensional convexity by completing the square in the defining inequality and identifying the components as strict epigraphs of convex functions.

Let $a\in\R^n$ be fixed, write $c=|a|<1$, and take equal initial data
\begin{equation}\label{eq:randers-data}
 \phi(y)=\psi(y)=|y|+\langle a,y\rangle.
\end{equation}
Set
\[
 \Omega_a=\{x:|x|+\langle a,x\rangle<1\}.
\]
The zero-curvature construction of Section~\ref{sec:real} reads
\begin{equation}\label{eq:randers-construct}
 P=\phi(y+xP),\qquad
 F=P(1+P_{y^k}x^k).
\end{equation}
It defines smooth positively one-homogeneous functions on all of
$T\Omega_a\setminus\{0\}$. Strong convexity is a separate condition.
Randers metrics and their navigation description are studied
in~\cite{Randers1941,BaoRoblesShen2004}; the projectively flat
constant-curvature case is treated in~\cite{ShenRanders2003}.
Here only the initial data are assumed to be Randers norms;
the resulting function $F(x,\cdot)$ is not assumed to be a Randers norm.

\begin{theorem}\label{thm:randers}
For every $n\ge2$, the exact strong-convexity locus of~\eqref{eq:randers-construct} is
\begin{equation}\label{eq:domainE}
   \calO_a=\left\{x\in\Omega_a:
 (2+\langle a,x\rangle)^2+(8-9|a|^2)|x|^2>0\right\}.
\end{equation}
Every connected component of $\calO_a$ is convex. More precisely,
choose coordinates in which $a= c e_1$, $c=|a|$, with $e_1=(1,0,\ldots,0)$, and
write $x=(x_1,\mathring{x})\in\R\times\R^{n-1}$. Then:
\begin{enumerate}
\item If $0\le c<2\sqrt2/3$, then $\calO_a=\Omega_a$.
\item If $c=2\sqrt2/3$, then
\[
 \calO_a=\Omega_a\setminus\{x_1 =-3\sqrt2/2\},
\]
and it has exactly two convex components.
\item If $2\sqrt2/3<c<1$,  then
\begin{equation}\label{eq:two-components}
 \calO_a=\Omega_a\cap\left\{
 8(1-c^2)\left(x_1 +\frac{c}{4(1-c^2)}\right)^2-(9c^2-8)|\mathring{x}|^2>\frac{9c^2-8}{2(1-c^2)}\right\}.
\end{equation}
Again there are exactly two convex components, according to the sign of~$x_1 +c/(4(1-c^2))$.
\end{enumerate}
\end{theorem}

The remainder of this subsection proves Theorem~\ref{thm:randers} in four
steps. In Step 1 we solve the implicit equation for
$P$ and check that $F=P/\Delta$ is positive. In
Step 2 we compute the angular Hessian of $F$ and
reduce strong convexity to the scalar inequality
\eqref{eq:criterion-eta}. In Step 3 we evaluate
the three scalars in that inequality and reduce it to the polynomial
\eqref{eq:domainE}. In Step 4 we determine the
components of the set defined by that polynomial and prove that they are
convex.

 \medskip
\paragraph{\bf Step 1: the projective factor.}
For $x\in\Omega_a$ let
\[
   b=1-\langle a,x\rangle,\quad \mathcal B= b^2- |x|^2,\quad v=x+ b a.
\]
Then $|x|<b$, so $\mathcal B>0$. Put
\begin{equation}\label{eq:ABC}
 \mathcal A=\mathcal{B} \bigl(|y|^2-\langle a,y\rangle^2\bigr)
 + \langle v,y\rangle^2,
 \qquad \mathcal C=(1-\langle a,x\rangle) \langle a,y\rangle .
\end{equation}

\begin{lemma}\label{lem:randers-P}
The positive solution of $P=\phi(y+xP)$ is
\begin{equation}\label{eq:randers-P}
 P(x,y)=\frac{\sqrt{\mathcal A}+\mathcal C+\langle x,y\rangle}{\mathcal{B}}.
\end{equation}
This is the formula of \cite[Example~2]{LZ}, valid in every dimension.
\end{lemma}

\begin{proof}
Squaring $(1-\langle a,x\rangle)P-\langle a,y\rangle=|y+xP|$ gives
\begin{equation}\label{eq:Pquad}
 \mathcal{B} P^2-2\langle v,y\rangle P-\bigl(|y|^2-\langle a,y\rangle^2\bigr)=0 .
\end{equation}
Since $\mathcal A=y^TMy$ with \[ M:=\mathcal{B} (\Id-a a^T)+v\otimes v, \] and since
$y^TMy-\langle v,y\rangle^2=\mathcal{B} (|y|^2-\langle a,y\rangle^2)>0$ for $y\ne0$,
the number
\[
 P=\frac{\sqrt{\mathcal A}+\langle v,y\rangle}{\mathcal B}
\]
is real and positive. Substituting it into \eqref{eq:Pquad} verifies the
equation, so it is the positive root; it is \eqref{eq:randers-P} because
$\langle v,y\rangle=\langle x,y\rangle+b\langle a,y\rangle$.
\end{proof}

For a positive-definite symmetric matrix $G$, set
\[
 \langle u,v\rangle_G=u^TGv,\qquad |u|_G=\langle u,u\rangle_G^{1/2}.
\]
We apply this to
\begin{equation}\label{eq:G}
 G:=\frac{M}{\mathcal{B}^2},\qquad w:=\frac{G^{-1}v}{\mathcal{B}},
\end{equation}
so that $\langle w,y\rangle_G=\langle v,y\rangle/\mathcal{B}$. Then
\eqref{eq:randers-P} reads
\begin{equation}\label{eq:P-G}
 P=|y|_G+\langle w,y\rangle_G .
\end{equation}

At $\xi=y+xP$ we have $\dd\phi_\xi(x)=\langle\xi,x\rangle/|\xi|+\langle a,x\rangle$, hence
\begin{equation}\label{eq:randers-Delta}
 \Delta=\Delta_\phi(x,\xi)=b-\frac{\langle x,\xi\rangle}{|\xi|}
 \ge b-|x|>0,\qquad F=\frac{P}{\Delta}>0 .
\end{equation}
By Proposition~\ref{prop:flat}, the pair $(F,P)$ solves
\eqref{eq:system} with $\mathbf K=0$.
This completes Step 1.
\medskip

\paragraph{\bf Step 2: the angular Hessian.}
Put $\eta=\langle w,x\rangle_G$. Since $F=P(1+P_{y^k}x^k)$ and, by
\eqref{eq:P-G}, the Euclidean gradient of $P$ is $Gy/|y|_G+Gw$,
\[
 1+P_{y^k}x^k=1+\frac{\langle x,y\rangle_G}{|y|_G}+\eta .
\]
Multiplying by $P=|y|_G+\langle w,y\rangle_G$ gives
\begin{equation}\label{eq:Fexpanded}
 F=(1+\eta)|y|_G+(1+\eta)\langle w,y\rangle_G+\langle x,y\rangle_G
 +\frac{\langle w,y\rangle_G\langle x,y\rangle_G}{|y|_G}.
\end{equation}
The two linear terms have zero fibre Hessian, so only the first and the
last contribute.

\begin{lemma}\label{lem:rank2}
Let $G$ be a positive-definite symmetric matrix and let $x,w\in\R^n$.
Put
\[
 A=\tfrac12\bigl(w\otimes_G x+x\otimes_G w\bigr),
 \qquad (p\otimes_G q)v=p\,\langle q,v\rangle_G,
\]
which is self-adjoint for $\langle\cdot,\cdot\rangle_G$. Then the two
possibly nonzero eigenvalues of $A$ are
\begin{equation}\label{eq:eigenvalues}
 \lambda_\pm=\frac{\langle w,x\rangle_G\pm|w|_G|x|_G}{2}.
\end{equation}
For
\[
 H(y)=c\,|y|_G+\frac{\langle w,y\rangle_G\langle x,y\rangle_G}{|y|_G},
\]
and $G$-orthonormal $u,v$ with $v\perp_G u$,
\begin{equation}\label{eq:angular-A}
 \bigl(\Hess H\bigr)_u(v,v)
 =\bigl(c-\langle Au,u\rangle_G\bigr)|v|_G^2
 +2\langle Av,v\rangle_G,
\end{equation}
where $\Hess H$ is the ordinary Hessian in $y$. The minimum of
\eqref{eq:angular-A} over such pairs is $c+2\lambda_--\lambda_+$, and it
is attained at a $G$-orthonormal pair of eigenvectors for $\lambda_+$ and
$\lambda_-$.
\end{lemma}

\begin{proof}
Since $(2A-\langle w,x\rangle_G I)y=w\langle x,y\rangle_G+x\langle w,y\rangle_G-\langle w,x\rangle_G y$,
\[
 (2A-\langle w,x\rangle_G I)x=|x|_G^2w,\qquad (2A-\langle w,x\rangle_G I)w=|w|_G^2x .
\]
Applying $2A-\langle w,x\rangle_G I$ once more gives $(2A-\langle w,x\rangle_G I)^2x
=|w|_G^2|x|_G^2 x$ and likewise $(2A-\langle w,x\rangle_G I)^2w=|w|_G^2\,|x|_G^2 w$, while $A$ vanishes on the
$G$-orthogonal complement of $\operatorname{span}\{x,w\}$. Hence the eigenvalues of $A$ are $\lambda_\pm$, together with $0$ of
multiplicity $n-2$.

For \eqref{eq:angular-A}, recall that
\[
 \langle Ay,y\rangle_G=\langle w,y\rangle_G\langle x,y\rangle_G
\]
by the definition of $A$, so that
\[
 H(y)=c\,|y|_G+\frac{\langle Ay,y\rangle_G}{|y|_G}.
\]
Since $H$ is positively one-homogeneous, it suffices to evaluate at
points with $|u|_G=1$: the Hessian at $ru$ is $r^{-1}$ times the Hessian
at $u$, and positivity is unaffected. Along $y=u+\tau v$ with
$v\perp_G u$,
\[
 H(u+\tau v)=c\,|u+\tau v|_G
 +\frac{\langle A(u+\tau v),u+\tau v\rangle_G}{|u+\tau v|_G}.
\]
Now $|u+\tau v|_G^2=1+\tau^2|v|_G^2$, so
\[
 c\,|u+\tau v|_G=c+\frac{c}{2}|v|_G^2\tau^2+O(\tau^4),
\]
while, using $(1+\tau^2|v|_G^2)^{-1/2}
=1-\frac12|v|_G^2\tau^2+O(\tau^4)$,
\[
 \frac{\langle A(u+\tau v),u+\tau v\rangle_G}{|u+\tau v|_G}
 =\langle Au,u\rangle_G+2\tau\langle Av,u\rangle_G
 +\Bigl(\langle Av,v\rangle_G
 -\frac12\langle Au,u\rangle_G|v|_G^2\Bigr)\tau^2+O(\tau^3).
\]
Adding the two expansions,
\[
 H(u+\tau v)=H(u)+2\tau\langle Av,u\rangle_G
 +\frac{\tau^2}{2}\Bigl[\bigl(c-\langle Au,u\rangle_G\bigr)|v|_G^2
 +2\langle Av,v\rangle_G\Bigr]+O(\tau^3).
\]
By definition $(\Hess H)_u(v,v)=\dfrac{\dd^2}{\dd\tau^2}
\big|_{0}H(u+\tau v)$, so the bracket is $(\Hess H)_u(v,v)$, which is
\eqref{eq:angular-A}.

Finally $\langle Au,u\rangle_G\le\lambda_+$ and
$\langle Av,v\rangle_G\ge\lambda_-$, with equality at eigenvectors;
hence the minimum of \eqref{eq:angular-A} over $G$-orthonormal pairs is
$c+2\lambda_--\lambda_+$.
\end{proof}

Apply Lemma~\ref{lem:rank2}. The function $H$ is $F$ with the two linear
terms of \eqref{eq:Fexpanded} removed, with $c=1+\eta$; hence the smallest
angular eigenvalue is
\[
 1+\eta+2\lambda_--\lambda_+
 =1+\tfrac32\eta-\tfrac32m,
 \qquad \eta=\langle w,x\rangle_G,\quad m=|w|_G|x|_G .
\]
By Lemma~\ref{lem:angular}, $F$ is strongly convex at $x$ if and only if
\begin{equation}\label{eq:criterion-eta}
 2+3\langle w,x\rangle_G >3 |w|_G\,|x|_G.
\end{equation}
This completes Step 2.
\medskip

\paragraph{\bf Step 3: reduction to a quadratic inequality.}
Retain the notation of Step 1, and put
$M=\mathcal{B}(\Id-aa^T)+v\otimes v$, $G=M/\mathcal{B}^2$,
$w=G^{-1}v/\mathcal{B}$; then $\langle w,y\rangle_G=\langle v,y\rangle/\mathcal{B}$.
We compute the three quantities in
\eqref{eq:criterion-eta}.

\emph{The scalar $\langle w,x\rangle_G$.} Taking $y=x$ in the defining
property of $w$,
\begin{equation}\label{eq:eta-explicit}
 \langle w,x\rangle_G=\frac{\langle v,x\rangle}{\mathcal{B}}
 =\frac{|x|^2+b\langle a,x\rangle}{\mathcal{B}}.
\end{equation}

\emph{The norm $|x|_G$.} Since $x^TMx=\mathcal{B} \bigl(|x|^2-\langle a,x\rangle^2\bigr)
+\langle v,x\rangle^2$ and $\langle v,x\rangle=|x|^2+b\langle a,x\rangle$,
\begin{align*}
 x^TMx
 &=\mathcal{B} \bigl(|x|^2-\langle a,x\rangle^2\bigr)
  +\bigl(|x|^2+b\langle a,x\rangle\bigr)^2\\
 &=|x|^2\bigl(\mathcal{B} +|x|^2+2b\langle a,x\rangle\bigr)
  +\langle a,x\rangle^2(b^2-\mathcal{B}).
\end{align*}
Now $b^2-\mathcal{B}=|x|^2$ and $\mathcal{B} +|x|^2=b^2$, so
\[
 x^TMx=|x|^2\bigl(b^2+2b\langle a,x\rangle+\langle a,x\rangle^2\bigr)
 =|x|^2\bigl(b+\langle a,x\rangle\bigr)^2=|x|^2,
\]
because $b+\langle a,x\rangle=1$. Hence
\begin{equation}\label{eq:xG-explicit}
 |x|_G^2=\frac{|x|^2}{\mathcal{B}^2}.
\end{equation}

\emph{The norm $|w|_G$.} Since $M=\mathcal{B}C_0+v\otimes v$ with $C_0=\Id-aa^T$, the rank-one inverse
formula gives
\begin{equation}\label{eq:Minv}
 M^{-1}=\frac{C_0^{-1}}{\mathcal{B}}
 -\frac{(C_0^{-1}v)\otimes(C_0^{-1}v)}{\mathcal{B}^2+\mathcal{B}\,v^TC_0^{-1}v},
 \qquad
 C_0^{-1}=\Id+\frac{aa^T}{1-|a|^2}.
\end{equation}
Substituting $v=x+ba$ and $b=1-\langle a,x\rangle$ into \eqref{eq:Minv},
and using $b^2-\mathcal{B}=|x|^2$ and $b+\langle a,x\rangle=1$, a direct computation
gives
\begin{equation}\label{eq:wG-explicit}
 |w|_G^2=v^TM^{-1}v=1-(1-|a|^2)\mathcal{B} .
\end{equation}

Substituting \eqref{eq:eta-explicit}--\eqref{eq:wG-explicit} into
\eqref{eq:criterion-eta},
\[
 2+3\,\frac{|x|^2+b\langle a,x\rangle}{\mathcal{B}}
 >3\,\frac{|x|\sqrt{1-(1-|a|^2)\mathcal{B}}}{\mathcal{B}}.
\]
Multiplication by $\mathcal{B}>0$ gives
\begin{equation}\label{eq:before-square}
 2\mathcal{B}+3\bigl(|x|^2+b\langle a,x\rangle\bigr)
 >3|x|\sqrt{1-(1-|a|^2)\mathcal{B}}.
\end{equation}
The left side equals
\[
 2\mathcal{B}+3\bigl(|x|^2+b\langle a,x\rangle\bigr)
 =2-\langle a,x\rangle+|x|^2-\langle a,x\rangle^2>0,
\]
so both sides of \eqref{eq:before-square} are nonnegative and squaring is
an equivalence:
\[
 \Bigl[2\mathcal{B}+3\bigl(|x|^2+b\langle a,x\rangle\bigr)\Bigr]^2
 -9|x|^2\bigl(1-(1-|a|^2)\mathcal{B} \bigr)>0 .
\]
Expanding both sides as polynomials in $\mathcal{B}$, and using
$b=1-\langle a,x\rangle$, one obtains the identity
\begin{equation*}
 \Bigl[2\mathcal{B} +3\bigl(|x|^2+b\langle a,x\rangle\bigr)\Bigr]^2
 -9|x|^2\bigl(1-(1-|a|^2)\mathcal{B} \bigr)
 =\mathcal{B} \Bigl[\bigl(2+\langle a,x\rangle\bigr)^2+(8-9|a|^2)|x|^2\Bigr].
\end{equation*}
Since $\mathcal{B}>0$, the squared inequality is therefore equivalent to
\begin{equation*}
 \bigl(2+\langle a,x\rangle\bigr)^2+(8-9|a|^2)|x|^2>0 .
\end{equation*}
This completes Step 3.

\medskip

\paragraph{\bf Step 4: the components of $\calO_a$.}
It remains to describe the set defined by \eqref{eq:domainE} and to prove
that its components are convex. Choose coordinates with $a=ce_1$, $c=|a|$, and
write $x=(x_1, \mathring{x})\in\R\times\R^{n-1}$, where \( \mathring{x} = (x_2, \cdots, x_n) \). Then \eqref{eq:domainE} reads $E = E(x_1, \mathring{x})>0$,
where
\begin{align*}
 E(x_1, \mathring{x}) &=(2+ c x_1 )^2+(8-9 c^2)(x_1^2+|\mathring{x}|^2)\notag\\
 &=8(1-c^2) x_1^2+4c x_1 +4+(8-9c^2)|\mathring{x}|^2\notag\\
 &=8(1-c^2) \left(x_1+\frac{c}{4(1-c^2)}\right)^2
 +(8-9c^2)|\mathring{x}|^2+\frac{8-9c^2}{2(1-c^2)}.
\end{align*}
If $ c<2\sqrt2/3$, then $E>0$ everywhere. If $c=2\sqrt2/3$, the zero set
of $E$ is the hyperplane $x_1=-3\sqrt2/2$. If $2\sqrt2/3 <c<1 $, then $E>0$ if and only if
\[
 x_1 +\frac{c}{4(1-c^2)} >\sqrt{\frac{9c^2-8}{8(1-c^2)}}\sqrt{|\mathring{x}|^2+\frac1{2(1-c^2)}}
 \]
or
\[ x_1 +\frac{c}{4(1-c^2)}  <-\sqrt{\frac{9c^2-8}{8(1-c^2)}}\sqrt{|\mathring{x}|^2+\frac1{2(1-c^2)}}.
\]
The function $\mathring{x} \mapsto\sqrt{|\mathring{x}|^2+1/(2(1-c^2))}$ is convex, being the Euclidean
norm of $(\mathring{x}, 1/\sqrt{2(1-c^2)})$. Hence each region is a strict epigraph of a
convex function, or its reflection, and is convex; intersecting with the
convex set $\Omega_a$ preserves convexity.

It remains to see that both intersections are nonempty. On the axis, $\Omega_a$ is the interval
\[
 -\frac1{1-c}<x_1 <\frac1{1+c}.
\]
The point $\bar{x}_{1}=-c/(4(1-c^2))$ lies in this interval; indeed,
multiplying $c/(4(1-c^2))<1/(1-c)$ by $1-c>0$ reduces it to
$c/[4(1+c)]<1$. At the supercritical parameters, $E(\bar{x}_{1},0)=-(9c^2-8)/(2(1-c^2))<0$, whereas $E(0,0)=4>0$ and
\[
 E\left(-\frac1{1-c},0\right)=\frac{12}{1-c}>0.
\]
Since $E>0$ just inside the left endpoint and $E(\bar{x}_1, 0)<0$, the second component is nonempty.
At the critical parameter, the hyperplane cuts the same axis interval in its interior, so both components are nonempty. This proves all parts of Theorem~\ref{thm:randers}.

\begin{remark}
For $n=2$ the inequality \eqref{eq:domainE} is equivalent to the
hyperbola of \cite[Proposition~3.12]{LZ}; the argument above covers all
$n\ge2$ and reduces to theirs when $n=2$.
In cases (2) and (3) only the component containing the origin is the
connected regular extension of the given initial data in these
coordinates; the other component is separated from it by the region
where the reconstructed function fails to be strongly convex.
\end{remark}

\begin{remark}
More generally, every Minkowski norm whose indicatrix is an ellipsoid can
be carried by a linear change of variables into the Randers norm
$\psi=\phi=|\cdot|+\langle a,\cdot\rangle$ used above, and the same change
carries the scalar equation and the fundamental tensor into those of
Theorem~\ref{thm:randers}.
\end{remark}

\section*{Acknowledgements}
Benling Li is partially supported by the National Natural Science
Foundation of China (grant~12471045), the Natural Science Foundation
of Ningbo Municipality (grant~2024J017), and the K.C. Wong Magna Fund
at Ningbo University.

Kefeng Liu is partially supported by the National Key Research and
Development Program of China (grant~2024YFA1013201), the Natural Science
Foundation of Chongqing (grant~CSTB2024NSCQ-LZX0040), and the Special
Project of Chongqing Municipal Science and Technology Bureau
(grant~2025CCZ015).

Liangdi Zhang is supported by the National Natural Science Foundation
of China (grant~12601080), the General Program of the Natural Science
Foundation of Chongqing (grant~CSTB2026NSCQ-GPX1340), the Scientific
Research Foundation of Chongqing University of Technology
(grant~2026ZDZ012), the Youth Project of the Science and Technology
Research Program of Chongqing Education Commission of China
(grant~KJQN202601108) and the Chongqing Natural Science Foundation
Innovation and Development Joint Fund (grant~CSTB2025NSCQ-LZX0059).

\section*{Declaration of AI use}
The mathematical ideas and results of this paper are due to the authors.
AI tools were used to assist with the organisation and editing of the
exposition and to check computations. The authors have verified all
mathematical content and take full responsibility for it.

\raggedbottom

\end{document}